\documentclass[11pt]{amsart}

\usepackage{amscd}
\usepackage{amsmath}
\usepackage{graphicx}
\usepackage{amsfonts}
\usepackage{amssymb}
\begin{document}

\title[Multiplicative Lie derivations]
{Multiplicative Lie derivation of  triangular 3-matrix rings}
\author{Zhenhui Chen,  Jinchuan Hou }
\address{School of Mathematics, Taiyuan University of Technology, Taiyuan, 030024, P. R.
China} \email[Jinchuan Hou]{jinchuanhou@aliyun.com }\email[Zhenhui Chen]{2691324902@qq.com}
\begin{abstract}

A map $\phi$ on an associative ring is called a multiplicative Lie
derivation if $\phi([x,y])=[\phi(x),y]+[x,\phi(y)]$ holds for any
elements $x,y$, where $[x,y]=xy-yx$ is the Lie product. In the
paper, we discuss the multiplicative Lie derivations on the
triangular 3-matrix rings $\mathcal T={\mathcal T}_3(\mathcal R_i;
\mathcal M_{ij})$.  Under the standard assumption $Q_i\mathcal
Z(\mathcal T)Q_i=\mathcal Z(Q_i\mathcal T Q_i)$, $i=1,2,3$, we show
that every multiplicative Lie derivation $\varphi:\mathcal
T\to\mathcal T$ has the standard form $\varphi=\delta+\gamma$ with
$\delta $ a derivation and $\gamma$ a center valued map vanishing
each commutator.

\end{abstract}


\thanks{This work is partially supported by
National Natural Science Foundation of China (11671294)} \maketitle

\section{Introduction}

 Let $R$ be an associative ring. Let $[x,y]=xy-yx$ denote the Lie product of $x,y\in R$. An additive map $\varphi$ of $R$ is
 called a derivation if $\varphi(xy)=\varphi(x)y+x\varphi(y)$ for all $x, y\in R$; $\varphi: R\rightarrow R$ is called a Lie derivation if
$\varphi([x,y])=[\varphi(x),y]+[x,\varphi(y)]$ for all $x, y\in R$.
Clearly, each derivation is a Lie derivation. The Lie derivations on
various kinds of rings and Algebras have been studied intensively
(see \cite{WS,FL,KM,MM,FW,AR,PS,WSC}).  Cheung
 in \cite{WS} described the form of Lie derivation on the triangular algebra. In \cite{FW}, Lu and Jing proved that
 if $\delta: B(X)\rightarrow B(X)$ is a linear map satisfying $\delta([A,B])=[\delta(A),B]+[A,\delta(B)]$ for any
 $A,B\in B(X)$ with $AB=0$, then $\delta$ can be decomposed as $d+\tau$, where $d$ is a derivation of $B(X)$ and
 $\tau:B(X)\rightarrow \mathbb{C}I$ is a linear map vanishing at commutators $[A,B]$ with $AB=0$. In paper \cite{YY},
 Du and Wang gave a description of Lie derivations of generalized matrix algebras. In \cite{QH},  Qi and  Hou characterized
 Lie derivations on a von Neumann algebra $M$ without central summands of type $I_{1}$. More generally, it was shown in \cite{QH}
 that, for any scalar $\xi$,  a
  additive map $L$ on $M$ satisfies $L(AB-\xi BA)=L(A)B-\xi BL(A)+AL(B)-\xi L(B)A$ whenever $AB=0$
  if and only if
 there exists an additive derivation $\varphi$ such that, (1) if $\xi=1$, then $L=\varphi+f$, where $f$ is an
 additive  map from $M$ into its center vanishing on $[A,B]$ with $AB=0$; (2) if $\xi=0$, then $L(I)\in Z(M)$
 and $L(A)=\varphi(A)+L(I)A$ for all $A$; (3) if $\xi$ is rational and $\xi\neq0,1$, then $L=\varphi$; (4) if $\xi$
 is not rational, then $\varphi(\xi I)=\xi L(I)$ and $L(A)=\varphi(A)+L(I)A$.   Benkovi$\check{c}$ in \cite{DB}
 proved that, under certain conditions, each generalized Lie derivation of a triangular algebra $\mathcal{A}$ is the sum
 of a generalized Lie derivation and a central map which vanishes on all commutators of
 $\mathcal{A}$.

 In  some situation, the assumption of the linearity or additivity on maps may be omitted.
 Thus
 the notion of Lie derivations was generalized to that of  multiplicative Lie derivations, and more generally, multiplicative Lie $n$-derivations.
Let us recall the following sequence
 of polynomials: $P_{1}(x)=x$ and $P_{n}(x_{1},x_{2},\cdot\cdot\cdot,x_{n})=[P_{n-1}(x_{1},x_{2},\cdot\cdot\cdot,x_{n-1}),x_{n}]$
 for all integers $n\geq2$. Thus, $P_{2}(x_{1},x_{2})=[x_{1},x_{2}],
 P_{3}(x_{1},x_{2},x_{3})=[[x_{1},x_{2}],x_{3}]$, $\ldots$.
  Let $n\geq 2$ be an integer. A map (no additivity is assumed) $\varphi: R\rightarrow R$ is called a multiplicative Lie $n$-derivation if
  $$\varphi(P_{n}(x_{1},x_{2},\cdot\cdot\cdot,x_{n}))=\sum^{n}_{i=1}P_{n}(x_{1},\cdot\cdot\cdot,x_{i-1},\varphi(x_{i}),x_{i+1},\cdot\cdot\cdot,x_{n})$$
 holds for all $x_1,x_2,\ldots, x_n\in\mathcal R$.
Particularly, a multiplicative Lie 2-derivation is called a
multiplicative Lie derivation.

Let $\mathcal A, \mathcal B$ be unital rings (res. algebras over a
commutative unital ring $\mathcal R$) and $\mathcal M$ be a
$(\mathcal A, \mathcal B)$-bimodule, which is faithful as a left
$\mathcal A$-module and also as a right $\mathcal B$-module; that
is, for any $A\in {\mathcal A}$ and $B\in \mathcal B$, $A{\mathcal
M}={\mathcal M}B=\{0\}$ imply $A=0$ and $B=0$. Recall that the
associative ring (res. algebra over $\mathcal R$)
$${\mathcal U}=\mbox{\rm Tri}(\mathcal A, \mathcal M, \mathcal B)=\{
\left(\begin{array}{cc}
A & M\\
  0 & B
\end{array}\right): A\in \mathcal A, M\in \mathcal M, B\in \mathcal B \}$$
under the usual matrix operations will be called a triangular ring
(res. algebra). In \cite{WYY}, Yu and Zhang discussed the
multiplicative Lie derivations on a triangular algebra ${\mathcal
U}={\rm Tri}(\mathcal A,\mathcal M,\mathcal B)$ over a commutative
ring. Under a standard assumption   $$\pi_{\mathcal A}({\mathcal
Z}(\mathcal U))={\mathcal Z}(\mathcal A)\ {\rm and}\ \pi_{\mathcal
B}({\mathcal Z}(\mathcal U))={\mathcal Z}(\mathcal B),\eqno(1.1)$$
they showed that every multiplicative Lie derivation $\varphi
:\mathcal U\to\mathcal U$ has the standard form
$\varphi=\delta+\gamma$, where $\delta $ is a derivation and $\gamma
$ is a center valuated map sending each commutator to zero. Here
$\pi_{\mathcal A}$ and $\pi_{\mathcal B}$ are projections defined
respectively by $\pi_{\mathcal A}(\left(\begin{array}{cc}
A & M\\
  0 & B
\end{array}\right))=A$ and  $\pi_{\mathcal B}(\left(\begin{array}{cc}
A & M\\
  0 & B
\end{array}\right))=B$.

More generally, under the assumptions that    $\mathcal{U}={\rm
Tri}(A,M,B)$ is $(n-1)$-torsion free   satisfying Eq.(1.1), and, in
addition, for the case $n\geq 3$, ${\mathcal A}$ or ${\mathcal B}$
has the property that $[a,\mathcal R]\in{\mathcal Z}(\mathcal
R)\Rightarrow a\in{\mathcal Z}(\mathcal R)$,
 Benkovi\v{c} and Eremita   \cite{DD}   proved that every multiplicative Lie $n$-derivation $\varphi:\mathcal{T}\rightarrow\mathcal{U}$ has the standard
  form $\varphi=\delta+\gamma$. In \cite{YYW}, Yao Wang and Yu Wang discussed the same problem for multiplicative Lie $n$-derivations on a certain class of generalized matrix algebras.

As a generalization of the notion of triangular rings, for each
positive integer $k\geq 2$, in \cite{BF}, Ferreira introduced a new
class of rings  which called triangular $k$-matrix rings. The notion
of triangular rings coincides with the notion of triangular
$2$-matrix rings. The works of \cite{WYY,DD, YYW} motivate us to
discuss
 the problem: {\it for $m\geq 3$,  under what conditions,  every multiplicative Lie $n$-derivation on a triangular $m$-matrix rings will have the standard
 form?}
We are not able to solve the problem for general $m$ in the present
paper. As a start, we show that, under a standard assumption, every
multiplicative
   Lie derivation on a triangular $3$-matrix ring has the standard
   form. But our method used in this paper is not valid for the case
   $m\geq 4$.

The paper is organized as follows. In section 2 we recall the notion
of triangular 3-matrix rings, provide a kind of the triangular
3-matrix rings that are not triangular rings. The main result,
including a corollary, is also given in this section. Section 3 is
devoted to discussing the general structural properties of the
multiplicative Lie derivations. We show that, under certain more
relaxed assumption on the triangular 3-matrix ring, every
multiplicative Lie derivation $\varphi$ has the form
$\varphi=\delta+\gamma+\xi$, where $\delta$ is a derivation,
$\gamma$ is a center valued map vanishing each commutator and $\xi$
is a multiplicative Lie derivation with a very small range. Based on
the results in Section 3, we give our proof of the main result in
Section 4.

 \section{Main result}

Though, for each $k\geq 2$, triangular $k$-matrix rings was
introduced in \cite{BF}, we only recall the definition of triangular
3-matrix rings here.

Let $\mathcal R_{1},\mathcal R_{2},\mathcal R_{3}$ be unital rings
and $\mathcal M_{ij}$ be $(\mathcal R_{i},\mathcal R_{j})$-bimodules
with $\mathcal M_{ii}=\mathcal R_{i}$ for all $1\leq i\leq j\leq3$.
Let $\varphi_{ijk}:\mathcal M_{ij\bigotimes \mathcal R_{j}}\mathcal
M_{jk}\rightarrow \mathcal M_{ik}$ be $(\mathcal R_{i},\mathcal
R_{k})$-bimodules homomorphisms with $\varphi_{iij}:\mathcal
R_{i}\otimes_{ \mathcal R_{i}}\mathcal M_{ij}\rightarrow \mathcal
M_{ij}$ and $\varphi_{ijj}:\mathcal M_{ij}\otimes _{\mathcal
R_{j}}R_{j}\rightarrow \mathcal M_{ij}$ the canonical multiplicative
maps for all $1\leq i\leq j\leq3$. Write $ab =\varphi_{ijk}(a\otimes
b)$ for $a\in \mathcal M_{ij},b\in \mathcal M_{jk}$. Assume that
 $\mathcal M_{ij}$ is faithful as a left $\mathcal
R_{i}$-module and faithful as a right $R_{j}$-module for all $1\leq
i<j\leq3.$ Let $\mathcal{T}=\mathcal T_3 (\mathcal R_i; \mathcal
M_{ij})$ be the set
$$\mathcal{T}=\mathcal T_3 (\mathcal R_i; \mathcal M_{ij})=\{\left(
\begin{array}{ccc}
 r_{11}&{m_{12}}&{m_{13}}\\
{0}&r_{22}&{m_{23}}\\
{0}&{0}&r_{33} \end{array} \right):r_{ii}\in R_{i}, m_{ij}\in
M_{ij},1\leq i<j\leq3\}. \eqno(2.1)$$ Clearly, with the matrix
operations of addition and multiplication, $\mathcal{T}$ is a ring
 which
 is called a triangular $3$-matrix ring. Note that, for $m>3$, the
 additional assumption is needed so that ${\mathcal T}_m({\mathcal
 R}_i: {\mathcal M}_{i,j})$ to be a ring:   $a(bc)=(ab)c$ for all $a\in \mathcal M_{ik}, b\in \mathcal M_{kl}$
and $c\in \mathcal M_{lj}$ with $1\leq i\leq k\leq l\leq j\leq m$
(Ref. \cite{BF}).

Obviously, setting
$$\mathcal{T}_{ij}=\{(m_{kt}):m_{kt}=\left\{\begin{array}{lll}
m_{ij}, & \mbox{\rm if}& (k,t)=(i,j) \\
0       & \mbox{\rm if}& (k,t)\neq(i,j)
\end{array}\right., \ 1\leq i\leq j\leq 3\}\subset \mathcal T,$$   we can write $\mathcal{T}=\bigoplus_{1\leq i\leq j\leq3}\mathcal{T}_{ij}$.
We also can identify $\mathcal T_{ij}$ with $\mathcal M_{ij}$.
Obviously, for any element $a_{ij}\in \mathcal{T}_{ij}$, we have
$a_{ij}a_{kj}=0$ whenever $j\neq k$.

Some canonical examples of triangular 3-matrix ring are upper
triangular matrix  rings ${\mathcal T}_n(\mathcal R) $ with $n\geq
3$ over a unital associative ring $\mathcal R$ and the nest algebras
with the nest containing more than 2 nontrivial elements.

It is easily seen that a triangular ring may not be a triangular
3-matrix ring. Conversely, the following example   shows that a
triangular 3-matrix ring may not be a triangular ring.

{\bf Example 2.1}  Let $\mathcal R$ be a unital associative ring
with the unit 1 and $M_{6}(\mathcal R)$ be the ring of all $6 \times
6$ matrices over $\mathcal R$. Let $\mathcal{A}\subseteq
M_{6}(\mathcal R)$ be the subset
$$\mathcal{A}=\{\left(
\begin{array}{cccccc}
 a_{11}&{a_{12}}&{0}&a_{14}&0&a_{16}\\
 a_{21}&{a_{22}}&{0}&a_{24}&0&a_{26}\\
 0&{0}&a&0&0&a_{36}\\
0&{0}&0&a&0&0\\
0&{0}&0&0&b&0\\
0&{0}&0&0&0&b\\\end{array} \right):a,b,a_{ij}\in \mathcal R\}.$$ It
is easily checked that $\mathcal{A}$ is a  subring of
$M_{6}(\mathcal R)$.

We claim that $\mathcal A$ is a triangular 3-matrix ring. To see this, let $ \mathcal A_1=M_2(\mathcal R)$, $\mathcal A_2=\mathcal A_3=\{\left(\begin{array}{cc} a &0\\
0 &a
\end{array}\right): a \in\mathcal R\}$, $\mathcal M_{12}=\mathcal M_{13}=\{\left(\begin{array}{cc} 0 &c\\ 0
&d
\end{array}\right): c,d\in \mathcal R\}$ and $\mathcal M_{23}=\{\left(\begin{array}{cc} 0 &c\\ 0
&0
\end{array}\right): c\in \mathcal R\}$; then
$ {\mathcal A}={\mathcal T}_3(\mathcal A_i; \mathcal M_{ij}).$ In
fact, it is obvious that $\mathcal A_i$s are unital rings and
$\mathcal M_{ij}$s are $(\mathcal A_i,\mathcal A_j)$-bimodules. We
check that $\mathcal M_{ij}$ is faithful as a left $\mathcal
A_i$-module and as a right $\mathcal A_j$-module. As $\mathcal
A_2=\mathcal A_3$ and $\mathcal M_{12}=\mathcal M_{13}$, we need
only to check this for $(i,j)=(1,2)$ and $(i,j)=(2,3)$.

Given any $A_1=\left(
\begin{array}{cc}
 a_{11}&{a_{12}}\\
a_{21}&a_{22}\end{array} \right)\in \mathcal A_{1}$. If
$A_1M_{12}=\left(
\begin{array}{cc}
 a_{11}&{a_{12}}\\
a_{21}&a_{22}\end{array} \right)\left(
\begin{array}{cc}
 0&c\\
0&d\end{array} \right)=0 $ holds for every $M_{12}=\left(
\begin{array}{cc}
 0&c\\
0&d\end{array} \right)\in \mathcal M_{12}$, then   $a_{ij}=0$ holds
for all $1\leq i,j\leq 2$ by taking $(c,d)=(1,0)$ or $(0,1)$, which
gives $A_1=0$. So $\mathcal M_{12}$ is faithful as a left ${\mathcal
A}_1$-module. It is obvious that $\mathcal M_{12}$ is faithful as a
right $\mathcal A_2$-module, and $\mathcal M_{23}$ is faithful.
Hence, $\mathcal A$ is a triangular 3-matrix ring.

We assert that $\mathcal A$ is not a triangular ring.

Clearly, one can not regard $\mathcal A$ as a triangular ring ${\rm
Tri}(\mathcal B_1,\mathcal M,\mathcal B_2)$ by $(3 + 3)$-partition,
that is, by letting
$${\mathcal B}_1=\{ \left(\begin{array}{ccc} a_{11} &a_{12} & 0 \\
a_{21} & a_{22} &0 \\ 0 & 0 & a \end{array}\right) :
a,a_{ij}\in\mathcal R \}
$$
$$\mathcal M=\{ \left(\begin{array}{ccc} a_{14} &0 & a_{16} \\
a_{24} & 0 &a_{26} \\ 0 & 0 & a_{36} \end{array}\right) :
a_{ij}\in\mathcal R \}
$$
and
$${\mathcal B}_2=\{ \left(\begin{array}{ccc} a  &0 & 0 \\
0 & b &0 \\ 0 & 0 & b \end{array}\right) : a, b\in\mathcal R \}.
$$
In fact, $\mathcal A$ is a proper subring of ${\rm Tri}(\mathcal
B_1,\mathcal M,\mathcal B_2)$ which contains elements $x=(a_{ij}) $
so that $a_{33}\not=a_{44}$. For the same reason, $\mathcal A$ is
not a triangular ring by $(5 + 1)$-partition.

Thus there are two possible ways, that is, $(2+4)$-partition and
$(4+2)$-partition, that might make $\mathcal A$ into  a triangular
ring.

{\bf Way 1}. \ $\mathcal A={\rm Tri}({\mathcal B}, \mathcal M,
\mathcal A_3)$, where ${\mathcal B}={\rm Tri}(\mathcal A_1, \mathcal
M_{12}, \mathcal A_2)$ and $$\mathcal M=\{ \left(
\begin{array}{cc}
 0&b\\
0&c\\
0&d\\
0&0 \end{array} \right) : b,c,d\in\mathcal R\}.$$ However, it is
easily seen that $\mathcal M$ is not faithful as a left ${\mathcal
B}$-module   by taking $$B=\left(
\begin{array}{cccc}
 a_{11}&{a_{12}}&0& 1 \\
a_{21}&a_{22}&0&a_{24}\\
0&0&a&0\\
0&0&0&a\end{array} \right)\in\mathcal B.$$

{\bf Way 2}. \ $\mathcal A={\rm Tri}({\mathcal A}_1, \mathcal N,
\mathcal C)$, where $\mathcal C={\rm Tri}({\mathcal A}_2, \mathcal
M_{23}, \mathcal A_3)$ and
$$\mathcal N=\{ \left(
\begin{array}{cccc}
 0&b &
0&d\\
0&c& 0&e \end{array} \right) : b,c,d,e\in\mathcal R\}.$$ But, by
taking $$C=\left(
\begin{array}{cccc}
 a &0 &0& 1\\
0&a &0&0\\
0&0&a^\prime &0\\
0&0&0&a^\prime \end{array} \right)\in \mathcal C$$  we see that
$\mathcal N$ is not faithful as a right ${\mathcal C}$-module.

Therefore, $\mathcal A$ is a triangular 3-matrix ring but not a
triangular ring. \hfill$\Box$

Before the statement of the main result, we introduce more
notations.

For any ring $\mathcal R$, $\mathcal Z (\mathcal R)$ stands for the
center of $\mathcal R$; that is, $\mathcal Z (\mathcal R)=\{Z :
Z\in\mathcal R\ {\rm and}\ [Z,\mathcal R]=0\}$. Let $\mathcal
T={\mathcal T}_3(\mathcal R_i;\mathcal M_{ij})$ be a triangular
3-matrix ring as defined in Eq.(2.1). There are three standard
idempotent elements in $\mathcal T$:
$$Q_{1}=\left(
\begin{array}{ccc}
 1&{0}&{0}\\
{0}&0&{0}\\
{0}&{0}&0 \end{array} \right),\quad Q_{2}=\left(
\begin{array}{ccc}
 0&{0}&{0}\\
{0}&1&{0}\\
{0}&{0}&0 \end{array} \right) \quad{\rm and}\quad Q_{3}=\left(
\begin{array}{ccc}
 0&{0}&{0}\\
{0}&0&{0}\\
{0}&{0}&1 \end{array} \right).$$ $\mathcal T$ is unital with the
unit $I=Q_1+Q_2+Q_3$. Denote by $\mathcal T_{ij}=Q_i{\mathcal
T}Q_j$, $i,j=1,2,3$. Clearly, ${\mathcal T}_{ij}=0$ if $i>j$,
${\mathcal T}_{i}={\mathcal T}_{ii}\cong {\mathcal R}_i$ and
${\mathcal T}_{ij}\cong {\mathcal M}_{ij}$ if $i<j$. Moreover,
$\mathcal T=\sum _{1\leq i\leq j\leq 3} {\mathcal T}_{ij}$. In this
paper, if  no confusion occurs, we identify $\mathcal T_i$ with
$\mathcal R_i$ and $\mathcal T_{ij}$ with $\mathcal M_{ij}$.

 {\bf Theorem 2.2.}    {\it Let $\mathcal{T}={\mathcal T}_3(\mathcal
R_i;\mathcal M_{ij})$ be a triangular 3-matrix ring. Assume that
$Q_{i}{\mathcal Z}(\mathcal{T})Q_{i}={\mathcal Z}(\mathcal{T}_{i})$,
$i=1,2,3$. Then each multiplicative Lie derivation $\varphi :
\mathcal{T}\rightarrow\mathcal{T}$ has the standard form
$$\varphi=\delta+\gamma,$$ where $\delta:
\mathcal{T}\rightarrow\mathcal{T}$ is a  derivation, $\gamma:
\mathcal{T}\rightarrow {\mathcal Z}(\mathcal{T})$ is a center valued
map such that $\gamma([\mathcal{T},\mathcal{T}])=0$.}

A proofs of Theorem 2.2   will be given in Section 4.

{\bf Corollary 2.3} \  {\it Let ${\mathcal A}={\rm T}_3({\mathcal
A}_i;{\mathcal M}_{ij})$ be the triangular 3-matrix ring constructed
from a ring $\mathcal R$ as in Example 2.1. Then every
multiplicative Lie derivation $\varphi :
\mathcal{A}\rightarrow\mathcal{A}$ has the  standard form. }

{\bf Proof.} It is obvious that ${\mathcal Z}(\mathcal
A_i)=\{\left(\begin{array}{cc} z &0 \\ 0 &z\end{array}\right) :
z\in{\mathcal Z}(\mathcal R)\}$ for $i=1,2,3$. It is also easily
checked that
$${\mathcal Z}(\mathcal A)=\{{\rm diag} (z,z,z,z,z, z) :
z\in{\mathcal Z}(\mathcal R)\}.
$$

As $Q_1={\rm diag}(1,1,0,0,0,0)$, $Q_2={\rm diag}(0,0,1,1,0,0)$ and
$Q_3={\rm diag}(0,0,0,0,1,1)$, we see that $Q_{i}{\mathcal
Z}(\mathcal{A})Q_{i}={\mathcal Z}(Q_i\mathcal{A}Q_{i})$, $i=1,2,3$.
Hence $\mathcal A$ satisfies the hypotheses of Theorem 2.2 and the
corollary follows. \hfill$\Box$

\section{Structure of multiplicative Lie derivations on triangular 3-matrix rings}

To prove Theorem 2.2, in this section, we consider the structure of
multiplicative Lie derivations on triangular 3-matrix rings in some
more general situation.

Denote by $\mathcal T$ the triangular 3-matrix ring ${\mathcal
T}_3(\mathcal R_i;\mathcal M_{ij})$. It is easily checked that the
center of $\mathcal{T}$ is the set
$$\mathcal Z(\mathcal{T})=\{{\rm diag}( r_{11}, r_{22},r_{33}):
r_{ii}\in{\mathcal R}_i,\ r_{ii}m_{ij}=m_{ij}r_{jj}\ \forall
m_{ij}\in \mathcal M_{ij},1\leq i<j\leq3\}. \eqno(3.1)$$
 Let $\mathcal{T}$ be the triangular 3-matrix
ring. It is clear from Eq.(3.1) that, for any $x\in\mathcal{T}$,
$$[x,\mathcal{T}]\subseteq \mathcal Z(\mathcal{T})\Rightarrow x\in
\mathcal Z(\mathcal{T}). $$

Denote by   $Q_i^\prime = I-Q_i$, $i=1,2,3$. Write $\mathcal
T_i=Q_i\mathcal T Q_i$, $\mathcal T_i^\prime=Q_i^\prime \mathcal T
Q_i^\prime$, $i=1,2,3$. Also write ${\mathcal M}_1=Q_1{\mathcal
T}Q_1'$, ${\mathcal M}_3=Q_3'{\mathcal T}Q_3$. Obviously, $\mathcal
T=\mathcal T_1+\mathcal M_1+\mathcal T_1'=\mathcal T_3'+\mathcal
M_3+\mathcal T_3$.

Denote by ${\mathcal V}_{23}=\{ w_{23}\in\mathcal M_{23} :
m_{12}w_{23}=0 \ \forall m_{12}\in{\mathcal M}_{12}\}$ and
${\mathcal V}_{12}=\{ w_{12}\in\mathcal M_{12} : w_{12}m_{23}=0 \
\forall m_{23}\in{\mathcal M}_{23}\}$. Also, denote by
$$
{\mathcal S}_{1}=\{ \left(\begin{array}{ccc} 0&0&0\\ 0&0& w_{23}\\
0&0&0 \end{array}\right) : w_{23}\in{\mathcal V}_{23}\}$$ and
$$ {\mathcal S}_{3}=\{ \left(\begin{array}{ccc} 0&w_{12}&0\\ 0&0& 0\\
0&0&0 \end{array}\right) : w_{12}\in{\mathcal V}_{12}\}.
$$

The following is our main result in this section which is needed to
prove Theorem 2.2.

{\bf Theorem 3.1.}    {\it Let $\mathcal{T}={\mathcal T}_3(\mathcal
R_i;\mathcal M_{ij})$ be a triangular 3-matrix ring and let
$k\in\{1,3\}$. Assume that $Q_{k}{\mathcal
Z}(\mathcal{T})Q_{k}={\mathcal Z}(\mathcal{T}_{k})$ and
$Q'_{k}{\mathcal Z}(\mathcal{T})Q'_{k}={\mathcal
Z}(\mathcal{T}'_{k})$. If every multiplicative Lie derivation from
${\mathcal T}_k'$ into itself has the standard form, then each
multiplicative Lie derivation $\varphi :
\mathcal{T}\rightarrow\mathcal{T}$ has the form
$$\varphi=\delta+\gamma+\xi,$$ where $\delta:
\mathcal{T}\rightarrow\mathcal{T}$ is a  derivation, $\gamma:
\mathcal{T}\rightarrow {\mathcal Z}(\mathcal{T})$ is a center valued
map such that $\gamma([\mathcal{T},\mathcal{T}])=0$ and $\xi :
\mathcal{T}\rightarrow\mathcal{S}_{k}$ is a multiplicative
 Lie derivation.}

We will give the details of the proof of Theorem 3.1 for the case
$k=1$. The case $k=3$ is treated with similarly.

\subsection{Proof of Theorem 3.1: Case $k=1$} In this case,
 $\mathcal{T}=\mathcal{T}_{1}+\mathcal M_{1}+\mathcal{T}'_{1}$ and we have $Q_{1}{\mathcal
Z}(\mathcal{T})Q_{1}={\mathcal Z}(\mathcal{T}_{1})$ and
$Q'_{1}{\mathcal Z}(\mathcal{T})Q'_{1}={\mathcal
Z}(\mathcal{T}'_{1})$.

 For an arbitrary element $x\in\mathcal T$,
one can then  write $x=a+m+b$ with $a\in\mathcal{T}_{1}, m\in
\mathcal M_{1}$ and $b\in\mathcal{T}'_{1}$. Let us define
projections
$\pi_{\mathcal{T}_{1}}:\mathcal{T}\rightarrow\mathcal{T}_{1}$ and
$\pi_{\mathcal{T}'_{1}}:\mathcal{T}\rightarrow\mathcal{T}'_{1}$ by
$\pi_{\mathcal{T}_{1}}(a+m+b)=a$ and
$\pi_{\mathcal{T}'_{1}}(a+m+b)=b$.  Obviously,
$\pi_{\mathcal{T}_{1}}(\mathcal Z(\mathcal{T}))\subseteq \mathcal
Z(\mathcal{T}_{1})$ and $\pi_{\mathcal{T}'_{1}}(\mathcal
Z(\mathcal{T}))\subseteq \mathcal Z(\mathcal{T}'_{1})$, and, there
exists a unique ring isomorphism
$\tau:\pi_{\mathcal{T}_{1}}(\mathcal Z(\mathcal{T}))\rightarrow
\pi_{\mathcal{T}'_{1}}(\mathcal Z(\mathcal{T}))$ such that
$am=m\tau(a)$ for all $m\in \mathcal M_{1}, a\in
\pi_{\mathcal{T}_{1}}(\mathcal Z(\mathcal{T}))$.

 In the sequel we assume that $\varphi : \mathcal
T\to\mathcal T$ is a  multiplicative Lie derivation. We prove
Theorem 3.1 for the case $k=1$ by a series lemmas.

The first lemma is an analogue of \cite[Lemma 3.1]{DD}.

{\bf Lemma 3.2.} {\it There exist an inner derivation
$d_{1}:\mathcal{T}-\mathcal{T}$ and a multiplicative Lie derivation
$\varphi_{1}:\mathcal{T}\rightarrow\mathcal{T}$ such that
$$\varphi=d_{1}+\varphi_{1} \quad {\rm and}\quad  Q_{1}\varphi_{1}(Q'_{1})Q'_{1}=0.$$}

{\bf Proof.} We define maps
$d_{1},\varphi_{1}:\mathcal{T}-\mathcal{T}$ by
$d_{1}(x)=[\varphi(Q'_{1}),x]$ and
$\varphi_{1}(x)=\varphi(x)-d_{1}(x)$. Obviously, $d_{1}$ is an inner
derivation and $\varphi_{1}$ is a multiplicative Lie derivation.
Since
$$\begin{array}{rl}\varphi_{1}(Q'_{1})&=\varphi(Q'_{1})-d_{1}(Q'_{1})\\
&=\varphi(Q'_{1})-[\varphi(Q'_{1}),Q'_{1}]\\
&=\varphi(Q'_{1})-\varphi(Q'_{1})Q'_{1}+Q'_{1}\varphi(Q'_{1}),
\end{array}$$
multiplying by $Q_{1} $ from the left and by $ Q'_{1}$ from the
right of the equation, we get
$Q_{1}\varphi_{1}(Q'_{1})Q'_{1}=0$.\hfill$\Box$

Thus, in Lemmas 3.3-3.5 we assume that $\varphi_1 :\mathcal T\to
\mathcal T$ is a multiplicative Lie derivations satisfying $
Q_{1}\varphi_{1}(Q'_{1})Q'_{1}=0$

 {\bf Lemma 3.3.} {\it For any $a\in
\mathcal{T}_{1}, b\in\mathcal{T}'_{1}$ and $ m\in\mathcal M_{1}$,
the following statements are true:}

(a)   $\varphi_{1}(Q_{1}),\ \varphi_{1}(Q'_{1})\in \mathcal
R_{1}\oplus\mathcal R_{2}\oplus\mathcal R_{3}$.

(b)  $Q_{1}\varphi_{1}(a)Q'_{1}=Q_{1}\varphi_{1}(b)Q'_{1}=0.$

(c)
$\varphi_{1}(a)=Q_{1}\varphi_{1}(a)Q_{1}+Q'_{1}\varphi_{1}(a)Q'_{1}$,
$\varphi_{1}(b)=Q_{1}\varphi_{1}(b)Q_{1}+Q'_{1}\varphi_{1}(b)Q'_{1}$
and $\varphi_{1}(m)=Q_{1}\varphi_{1}(m)Q'_{1}.$

(d)  $Q_{1}\varphi_{1}(\mathcal{T}'_{1})Q_{1}\subseteq \mathcal
Z(\mathcal{T}_{1})$ and $
Q'_{1}\varphi_{1}(\mathcal{T}_{1})Q'_{1}\subseteq \mathcal
Z(\mathcal{T}'_{1}).$

{\bf Proof.} It is obvious that $\varphi(0)=0$ and
$\varphi(I)\in\mathcal Z(\mathcal T)$.

(a) Since $Q_{1}\varphi_{1}(Q'_{1})Q'_{1}=0$,
$$\varphi_{1}(Q'_{1})\in\{\left(
\begin{array}{ccc}
 a_{11}&{0}&{0}\\
{0}&a_{22}&{a_{23}}\\
{0}&{0}&a_{33} \end{array} \right):a_{ii}\in \mathcal
R_{i}(i=1,2,3),a_{23}\in \mathcal M_{23}\}.$$ Note that
$$0=\varphi([Q_{1},Q'_{1}])=\varphi(Q_{1})Q'_{1}-Q'_{1}\varphi(Q_{1})+Q_{1}\varphi(Q'_{1})-\varphi(Q'_{1})Q_{1}.\eqno(3.2)$$
By Lemma 3.2, we have $$\begin{array}{rl}\varphi_{1}(Q_{1})&=\varphi(Q_{1})-d_{1}(Q_{1})\\
&=\varphi(Q_{1})-[\varphi(Q'_{1}),Q_{1}]               \\
&=\varphi(Q_{1})-(\varphi(Q'_{1})Q_{1}-Q_{1}\varphi(Q'_{1})).
\end{array}\eqno(3.3)$$ Thus Eq.(3.2) together with Eq.(3.3) gives
$$\varphi_{1}(Q_{1})=\varphi(Q_{1})-(\varphi(Q_{1})Q'_{1}-Q'_{1}\varphi(Q_{1})).$$
Writing $$\varphi(Q_{1})=\left(
\begin{array}{ccc}
 a_{11}&{a_{12}}&{a_{13}}\\
{0}&a_{22}&{a_{23}}\\
{0}&{0}&a_{33} \end{array} \right),$$ the above formula leads to
$$\begin{array}{rl} \varphi_{1}(Q_{1})  =&\left(
\begin{array}{ccc}
 a_{11}&{a_{12}}&{a_{13}}\\
{0}&a_{22}&{a_{23}}\\
{0}&{0}&a_{33} \end{array} \right)-\left(
\begin{array}{ccc}
 a_{11}&{a_{12}}&{a_{13}}\\
{0}&a_{22}&{a_{23}}\\
{0}&{0}&a_{33} \end{array} \right)\left(
\begin{array}{ccc}
 0&{0}&{0}\\
{0}&1&{0}\\
{0}&{0}&1 \end{array} \right) \\ &+\left(
\begin{array}{ccc}
 0&{0}&{0}\\
{0}&1&{0}\\
{0}&{0}&1 \end{array} \right)\left(
\begin{array}{ccc}
 a_{11}&{a_{12}}&{a_{13}}\\
{0}&a_{22}&{a_{23}}\\
{0}&{0}&a_{33} \end{array} \right)\\
=&\left(
\begin{array}{ccc}
 a_{11}&{a_{12}}&{a_{13}}\\
{0}&a_{22}&{a_{23}}\\
{0}&{0}&a_{33} \end{array} \right)-\left(
\begin{array}{ccc}
 0&{a_{12}}&{a_{13}}\\
{0}&a_{22}&{a_{23}}\\
{0}&{0}&a_{33} \end{array} \right)+\left(
\begin{array}{ccc}
 0&{0}&{0}\\
{0}&a_{22}&{a_{23}}\\
{0}&{0}&a_{33} \end{array} \right)               \\
=&\left(
\begin{array}{ccc}
 a_{11}&{0}&{0}\\
{0}&a_{22}&{a_{23}}\\
{0}&{0}&a_{33} \end{array} \right).                                \end{array}$$
Therefore $$\varphi_{1}(Q_{1})\in\{\left(
\begin{array}{ccc}
 a_{11}&{0}&{0}\\
{0}&a_{22}&{a_{23}}\\
{0}&{0}&a_{33} \end{array} \right):a_{ii}\in \mathcal
R_{i}(i=1,2,3), a_{23}\in \mathcal M_{23}\}.$$  Similarly one can
check that
$$\varphi_{1}(Q_{2})\in\{\left(
\begin{array}{ccc}
 a_{11}&{0}&{0}\\
{0}&a_{22}&{a_{23}}\\
{0}&{0}&a_{33} \end{array} \right):a_{ii}\in \mathcal
R_{i}(i=1,2,3),a_{23}\in \mathcal M_{23}\}.$$ So we can write
$$\varphi_{1}(Q_{1})=\left(
\begin{array}{ccc}
 a_{11}&{0}&{0}\\
{0}&a_{22}&{a_{23}}\\
{0}&{0}&a_{33} \end{array} \right),\quad \varphi_{1}(Q_{2})=\left(
\begin{array}{ccc}
 b_{11}&{0}&{0}\\
{0}&b_{22}&{b_{23}}\\
{0}&{0}&b_{33} \end{array} \right)$$ and
$$\varphi_{1}(Q'_{1})=\left(
\begin{array}{ccc}
 c_{11}&{0}&{0}\\
{0}&c_{22}&{c_{23}}\\
{0}&{0}&c_{33} \end{array} \right).$$
Since
$$\begin{array}{rl}0=&\varphi_{1}([Q_{1},Q_{2}])\\
=&\varphi_{1}(Q_{1})Q_{2}-Q_{2}\varphi_{1}(Q_{1})+Q_{1}\varphi_{1}(Q_{2})-\varphi_{1}(Q_{2})Q_{1}\\
=&\left(
\begin{array}{ccc}
 a_{11}&{0}&{0}\\
{0}&a_{22}&{a_{23}}\\
{0}&{0}&a_{33} \end{array} \right)\left(
\begin{array}{ccc}
 1&{0}&{0}\\
{0}&1&{0}\\
{0}&{0}&0 \end{array} \right)-\left(
\begin{array}{ccc}
 1&{0}&{0}\\
{0}&1&{0}\\
{0}&{0}&0 \end{array} \right)\left(
\begin{array}{ccc}
 a_{11}&{0}&{0}\\
{0}&a_{22}&{a_{23}}\\
{0}&{0}&a_{33} \end{array} \right)\\
&+\left(
\begin{array}{ccc}
 1&{0}&{0}\\
{0}&0&{0}\\
{0}&{0}&0 \end{array} \right)\left(
\begin{array}{ccc}
 b_{11}&{0}&{0}\\
{0}&b_{22}&{b_{23}}\\
{0}&{0}&b_{33} \end{array} \right)-\left(
\begin{array}{ccc}
 b_{11}&{0}&{0}\\
{0}&b_{22}&{b_{23}}\\
{0}&{0}&b_{33} \end{array} \right)\left(
\begin{array}{ccc}
 1&{0}&{0}\\
{0}&0&{0}\\
{0}&{0}&0 \end{array} \right)\\
=&\left(
\begin{array}{ccc}
 a_{11}&{0}&{0}\\
{0}&a_{22}&{0}\\
{0}&{0}&0 \end{array} \right) -\left(
\begin{array}{ccc}
 a_{11}&{0}&{0}\\
{0}&a_{22}&{a_{23}}\\
{0}&{0}&0 \end{array} \right) +\left(
\begin{array}{ccc}
 b_{11}&{0}&{0}\\
{0}&0&{0}\\
{0}&{0}&0 \end{array} \right)-\left(
\begin{array}{ccc}
 b_{11}&{0}&{0}\\
{0}&0&{0}\\
{0}&{0}&0 \end{array} \right)          \\
=&\left(
\begin{array}{ccc}
 0&{0}&{0}\\
{0}&0&{-a_{23}}\\
{0}&{0}&0 \end{array} \right),
\end{array}$$ we have $a_{23}=0$ and thus $$\varphi_{1}(Q_{1})=\left(
\begin{array}{ccc}
 a_{11}&{0}&{0}\\
{0}&a_{22}&{0}\\
{0}&{0}&a_{33} \end{array} \right)\in \mathcal R_{1}\oplus\mathcal
R_{2}\oplus\mathcal R_{3}.$$ Similarly
$$\begin{array}{rl}0=&\varphi_{1}([Q'_{1},Q_{2}])\\
=&\varphi_{1}(Q'_{1})Q_{2}-Q_{2}\varphi_{1}(Q'_{1})+Q'_{1}\varphi_{1}(Q_{2})-\varphi_{1}(Q_{2})Q'_{1}\\
=&\left(
\begin{array}{ccc}
 c_{11}&{0}&{0}\\
{0}&c_{22}&{c_{23}}\\
{0}&{0}&c_{33} \end{array} \right)\left(
\begin{array}{ccc}
 1&{0}&{0}\\
{0}&1&{0}\\
{0}&{0}&0 \end{array} \right)-\left(
\begin{array}{ccc}
 1&{0}&{0}\\
{0}&1&{0}\\
{0}&{0}&0 \end{array} \right)\left(
\begin{array}{ccc}
 c_{11}&{0}&{0}\\
{0}&c_{22}&{c_{23}}\\
{0}&{0}&c_{33} \end{array} \right)\\
&+\left(
\begin{array}{ccc}
 0&{0}&{0}\\
{0}&1&{0}\\
{0}&{0}&1 \end{array} \right)\left(
\begin{array}{ccc}
 b_{11}&{0}&{0}\\
{0}&b_{22}&{b_{23}}\\
{0}&{0}&b_{33} \end{array} \right)-\left(
\begin{array}{ccc}
 b_{11}&{0}&{0}\\
{0}&b_{22}&{b_{23}}\\
{0}&{0}&b_{33} \end{array} \right)\left(
\begin{array}{ccc}
 0&{0}&{0}\\
{0}&1&{0}\\
{0}&{0}&1 \end{array} \right)\\
=&\left(
\begin{array}{ccc}
 c_{11}&{0}&{0}\\
{0}&c_{22}&{0}\\
{0}&{0}&0 \end{array} \right) -\left(
\begin{array}{ccc}
 c_{11}&{0}&{0}\\
{0}&c_{22}&{c_{23}}\\
{0}&{0}&0 \end{array} \right) +\left(
\begin{array}{ccc}
 0&{0}&{0}\\
{0}&b_{22}&{b_{23}}\\
{0}&{0}&b_{33} \end{array} \right)-\left(
\begin{array}{ccc}
 0&{0}&{0}\\
{0}&b_{22}&{b_{23}}\\
{0}&{0}&b_{33} \end{array} \right)          \\
=&\left(
\begin{array}{ccc}
 0&{0}&{0}\\
{0}&0&{-c_{23}}\\
{0}&{0}&0 \end{array} \right)
\end{array}$$ entails $c_{23}=0$. Hence $$\varphi_{1}(Q'_{1})=\left(
\begin{array}{ccc}
 c_{11}&{0}&{0}\\
{0}&c_{22}&{0}\\
{0}&{0}&c_{33} \end{array} \right)\in \mathcal R_{1}\oplus\mathcal
R_{2}\oplus\mathcal R_{3},$$ too.

(b) Let $x\in \mathcal{T}_{1}\bigcup\mathcal{T}'_{1}$. Then $[x,
Q'_{1}]=Q_{1}xQ'_{1}=0$. Since $\varphi_{1}(0)=0$, we have
$$\begin{array}{rl}0&=\varphi_{1}(Q_{1}xQ'_{1})=\varphi_{1}([x, Q'_{1}])\\
&=[\varphi_{1}(x),Q'_{1}]+[x,\varphi_{1}(Q'_{1})]\\
&=Q_{1}\varphi_{1}(x)Q'_{1}+Q_{1}[x,\varphi_{1}(Q'_{1})] Q'_{1}.
\end{array} \eqno(3.4)$$ However the property  $Q_{1}\varphi_{1}(Q'_{1})Q'_{1}=0$
implies that
$$\begin{array}{rl}&Q_{1}[x,\varphi_{1}(Q'_{1})]Q'_{1}\\
=&Q_{1}(x\varphi_{1}(Q'_{1})-\varphi_{1}(Q'_{1})x)Q'_{1}            \\
=&Q_{1}(Q_{1}xQ_{1}\varphi_{1}(Q'_{1})-\varphi_{1}(Q'_{1})Q'_{1}xQ_{1})Q'_{1}               \\
=&0                                  \end{array}$$ holds for every
$x\in \mathcal T_1$. So
$Q_{1}[\mathcal{T}_{1},\varphi_{1}(Q'_{1})]Q'_{1}=0$. Similarly, we
have $Q_{1}[\mathcal{T}'_{1},\varphi_{1}(Q'_{1})]Q'_{1}=0$. So by
Eq.(3.4) we see that $Q_{1}\varphi_{1}(x)Q'_{1}=0$ for all $x\in
\mathcal{T}_{1}\bigcup\mathcal{T}'_{1}$; That is, (b) is true.

(c) By (b) it is clear that
$\varphi_{1}(a)=Q_{1}\varphi_{1}(a)Q_{1}+Q'_{1}\varphi_{1}(a)Q'_{1}
$ and
$\varphi_{1}(b)=Q_{1}\varphi_{1}(b)Q_{1}+Q'_{1}\varphi_{1}(b)Q'_{1}$
are true respectively for every $a\in{\mathcal T}_1$ and
$b\in{\mathcal T}_1^\prime$. Now, for each $m\in \mathcal M_{1}$, we
have
$$\begin{array}{rl}\varphi_{1}(m)&=\varphi_{1}[m,Q'_{1}]\\
&=[\varphi_{1}(m),Q'_{1}]+[m,\varphi_{1}(Q'_{1})]        \\
&=\varphi_{1}(m)Q'_{1}-Q'_{1}\varphi_{1}(m)+[m,\varphi_{1}(Q'_{1})].
\end{array}$$ Multiplying by $ Q_{1}$ from the left hand side and by $ Q'_{1}$
from the right hand side of the equation, we see that
$$[m,\varphi_{1}(Q'_{1})]=Q_1[m,\varphi_{1}(Q'_{1})]Q_1^\prime=0.$$
This implies that
$\varphi_{1}(m)=\varphi_{1}(m)Q'_{1}-Q'_{1}\varphi_{1}(m)=Q_{1}\varphi_{1}(m)Q'_{1}$.
So (c) is true.

 (d) Let $a\in\mathcal{T}_{1},b\in\mathcal{T}'_{1}$ and $ m\in\mathcal M_{1}$.
Since $[a,b]=0$, we have
$$[\varphi_{1}(a),b]+[a,\varphi_{1}(b)]=\varphi_{1}([a,b])=0.
$$ Applying (c), we see that
$$[a,Q_{1}\varphi_{1}(b)Q_{1}]=[a,\varphi_{1}(b)]=0$$ and
$$[Q'_{1}\varphi( a)Q'_{1},b]=[\varphi_{1}(a),b]=0$$
hold for all $a\in{\mathcal T}_1$ and all $b\in{\mathcal
T}_1^\prime$. Therefore, $
Q_{1}\varphi(\mathcal{T}'_{1})Q_{1}\subseteq Z(\mathcal{T}_{1})$ and
$Q'_{1}\varphi(\mathcal{T}_{1})Q'_{1}\subseteq Z(\mathcal{T}'_{1})$.
This completes the proof of (d).

Furthermore, we have

{\bf Lemma 3.4.} {\it  $\varphi_{1}(Q_{1}), \varphi_{1}(Q'_{1})\in
\mathcal Z(\mathcal{T})$ and
$\varphi_{1}(Q_{1}xQ'_{1})=Q_{1}\varphi_{1}(x)Q'_{1}$ for all
$x\in\mathcal T$.}

{\bf Proof.}
For any $x\in\mathcal{T}$, we have $$\begin{array}{rl}\varphi_{1}(Q_{1}xQ'_{1})&=\varphi_{1}([Q_{1},x])\\
&=[\varphi_{1}(Q_{1}),x]+[Q_{1},\varphi_{1}(x)].       \\
\end{array}\eqno(3.5)$$
 Replacing $x$ by $Q_{1}xQ'_{1}$ in (3.5), we  obtain
$$\varphi_{1}(Q_{1}xQ'_{1})=[\varphi_{1}(Q_{1}),Q_{1}xQ'_{1}]+[Q_{1},\varphi_{1}(Q_{1}xQ'_{1})].$$
By (c) one gets $[\varphi_{1}(Q_{1}),Q_{1}xQ'_{1}]=0$, which in turn
implies that
$$[\varphi_{1}(Q_{1}),\mathcal M_{1}]=0.\eqno(3.6)$$
Meanwhile, by (c), it is clear that
 \if false
$$\varphi_{1}(Q_{1}\mathcal{T}_{1}Q'_{1})=[\varphi_{1}(Q_{1}),\mathcal{T}_{1}]+[Q_{1},\varphi_{1}(\mathcal{T}_{1})]=0$$
and
$$\varphi_{1}(Q_{1}\mathcal{T}'_{1}Q'_{1})=[\varphi_{1}(Q_{1}),\mathcal{T}'_{1}]+[Q_{1},\varphi_{1}(\mathcal{T}'_{1})]=0,$$
according to (iii), which in turn implies
$$[Q_{1},\varphi_{1}(\mathcal{T}_{1})]=[Q_{1},Q_{1}\varphi_{1}(\mathcal{T}_{1})Q_{1}+Q'_{1}\varphi_{1}(\mathcal{T}_{1})Q'_{1}]=0$$
and
$$[Q_{1},\varphi_{1}(\mathcal{T}'_{1})]=[Q_{1},Q_{1}\varphi_{1}(\mathcal{T}'_{1})Q_{1}+Q'_{1}\varphi_{1}(\mathcal{T}'_{1})Q'_{1}]=0.$$
\fi
 $$[\varphi_{1}(Q_{1}),\mathcal{T}_{1}]=0,\quad{\rm and}\quad
[\varphi_{1}(Q_{1}),\mathcal{T}'_{1}]=0.\eqno(3.7)$$ Now (3.6) and
(3.7) together ensures that
$$\varphi_{1}(Q_{1})\in \mathcal Z(\mathcal{T}).$$
An analogous  argument shows that \if false , $$\begin{array}{rl}\varphi_{1}(Q_{1}xQ'_{1})&=\varphi_{1}([x,Q'_{1}])\\
&=[\varphi_{1}(x),Q'_{1}]+[x,\varphi_{1}(Q'_{1})]       \\
\end{array}\eqno(6)$$
for all $x\in\mathcal{T}$, Replacing $x$ by $Q_{1}xQ'_{1}$ in (6), we can get
$$\varphi_{1}(Q_{1}xQ'_{1})=[\varphi_{1}(Q_{1}xQ'_{1}),Q'_{1}]+[Q_{1}xQ'_{1},\varphi_{1}(Q'_{1})],$$
so $[Q_{1}xQ'_{1},\varphi_{1}(Q'_{1})]=0$, which in turn implies $$[M_{1},\varphi_{1}(Q'_{1})]=0.\eqno(7)$$
Meanwhile, $$\varphi_{1}(Q_{1}\mathcal{T}_{1}Q'_{1})=[\varphi_{1}(\mathcal{T}_{1}),Q'_{1}]+[\mathcal{T}_{1},\varphi_{1}(Q'_{1})]=0$$
and $$\varphi_{1}(Q_{1}\mathcal{T}'_{1}Q'_{1})=[\varphi_{1}(\mathcal{T}'_{1}),Q'_{1}]+[\mathcal{T}'_{1},\varphi_{1}(Q'_{1})]=0,$$
according to (iii), it implies
$$[\varphi_{1}(\mathcal{T}_{1}),Q'_{1}]=[Q_{1}\varphi_{1}(\mathcal{T}_{1})Q_{1}+Q'_{1}\varphi_{1}(\mathcal{T}_{1})Q'_{1},Q'_{1}]=0$$
and
$$[\varphi_{1}(\mathcal{T}'_{1}),Q'_{1}]=[Q_{1}\varphi_{1}(\mathcal{T}'_{1})Q_{1}+Q'_{1}\varphi_{1}(\mathcal{T}'_{1})Q'_{1},Q'_{1}]=0.$$
Thus
$$[\mathcal{T}_{1},\varphi_{1}(Q'_{1})]=0,[\mathcal{T}'_{1},\varphi_{1}(Q'_{1})]=0.\eqno(8)$$
Now, (7) together with (8), we obtain
$[\mathcal{T},\varphi_{1}(Q'_{1})]=0,$ so \fi
$$\varphi_{1}(Q'_{1})\in \mathcal Z(\mathcal{T}),$$ too.
Consequently, $$\begin{array}{rl}\varphi_{1}(Q_{1}xQ'_{1})&=\varphi_{1}([Q_{1},x])\\
&=[\varphi_{1}(Q_{1}),x]+[Q_{1},\varphi_{1}(x)]       \\
&=Q_{1}\varphi_{1}(x)Q'_{1}
\end{array}$$
holds for all $x\in\mathcal T$. This establishes the
lemma.\hfill$\Box$

{\bf Lemma 3.5.} {\it For any $a\in \mathcal{T}_{1}, b\in
\mathcal{T}'_{1}$ and $ m\in\mathcal M_{1}$, we have
$$\varphi_{1}(am)=\varphi_{1}(a)m+a\varphi_{1}(m)-m\varphi_{1}(a)$$ and $$ \varphi_{1}(mb)=\varphi_{1}(m)b+m\varphi_{1}(b)-\varphi_{1}(b)m.$$}

{\bf Proof.} For any $a\in \mathcal{T}_{1}, b\in \mathcal{T}'_{1}$
and $ m\in\mathcal M_{1}$, by applying Lemma 3.3, we have
$$\begin{array}{rl}\varphi_{1}(am)&=\varphi_{1}([a,m])\\
&=[\varphi_{1}(a),m]+[a,\varphi_{1}(m)]\\
&=\varphi_{1}(a)m-m\varphi_{1}(a)+a\varphi_{1}(m)-\varphi_{1}(m)a\\
&=\varphi_{1}(a)m+a\varphi_{1}(m)-m\varphi_{1}(a)  \end{array}$$
  and
$$\begin{array}{rl}\varphi_{1}(mb)&=\varphi_{1}([m,b])\\
&=[\varphi_{1}(m),b]+[m,\varphi_{1}(b)]\\
&=\varphi_{1}(m)b-b\varphi_{1}(m)+m\varphi_{1}(b)-\varphi_{1}(b)m\\
&=\varphi_{1}(m)b+m\varphi_{1}(b)-\varphi_{1}(b)m.  \end{array}$$
So, the lemma is true. \hfill $\Box$

{\bf Lemma 3.6.} {\it For any $a\in\mathcal{T}_{1}$, $
b\in\mathcal{T}'_{1}$ and $ m\in\mathcal M_{1}$, we have}

(1) {\it Both $\varphi_{1}(a+m)-\varphi_{1}(a)-\varphi_{1}(m)$ and
$\varphi_{1}(a+b)-\varphi_{1}(a)-\varphi_{1}(b)\in\mathcal
Z(\mathcal{T}) $.}

(2) {\it $\varphi_{1}$~is additive on $\mathcal M_{1}$}.

(3) {\it With ${\mathcal V}_{23}=\{ w_{23}\in{\mathcal M}_{23} :
{\mathcal M}_{12}w_{23}=0\}$, we have
$$\begin{array}{rl} &
\varphi_{1}(a+m+b)-\varphi_{1}(a)-\varphi_{1}(b)-\varphi_{1}(m) \\
\in & \{ {\small \left(
\begin{array}{ccc}
 a_{11}&0&0\\
{0}&a_{22}&w_{23}\\
{0}&{0}&a_{33} \end{array} \right)} : w_{23}\in{\mathcal V}_{23},
a_{11}\in{\mathcal Z}({\mathcal T}_{1}), a_{11}m_{1j}=m_{1j}a_{jj} \
\forall m_{1j}\in\mathcal M_{1j}, j=2,3\}.\end{array}$$ }

\textbf{Proof.} (1)  Let $a\in\mathcal{T}_{1}$ and $m\in\mathcal
M_{1}$. For any $m'\in\mathcal M_1$, since $[a+m,m']=[a,m']$, we
have
$$\varphi_{1}([a,m'])=\varphi_{1}([a+m,m']).$$
Due to the fact that $\varphi_{1}(m)\in\mathcal M_{1}$,
 we have $[\varphi_{1}(m),m']=0$, and then
$$[\varphi_{1}(a+m),m']=[\varphi_{1}(a),m'].$$
So
$$[\varphi_{1}(a+m)-\varphi_{1}(a)-\varphi_{1}(m),m']=0 \eqno(3.8)$$
holds for all   $m'  \in\mathcal M_1$. Notice that, for any
$b'\in{\mathcal T}_1^\prime$,
$$\begin{array}{rl}\varphi_{1}(mb')&=\varphi_{1}([a+m,b'])\\
&=[\varphi_{1}(a+m),b']+[a+m,\varphi_{1}(b')],
 \end{array}$$
and, by Lemma 3.5,
$$\begin{array}{rl}\varphi_{1}(mb')&=\varphi_{1}(m)b'+m\varphi_{1}(b')-\varphi_{1}(b')m\\
&=[\varphi_{1}(m),b']+[m,\varphi_{1}(b')].
 \end{array}$$
The above two equations, together with the facts
$[a,\varphi_{1}(b')]=0$ and $[\varphi_{1}(a),b']=0$ (see Lemma 3.3
(d)), entail that
$$[\varphi_{1}(a+m)-\varphi_{1}(m)-\varphi_{1}(a),b']=0\eqno(3.9)$$
holds for all $b'\in{\mathcal T}_1^\prime$.

Next, by Lemma 3.4, we have
$Q_{1}\varphi_{1}(a+m)Q'_{1}=\varphi_{1}(m)$. Thus we may write
$$\varphi_{1}(a+m)-\varphi_{1}(a)-\varphi_{1}(m)=\left(
\begin{array}{ccc}
 a_{11}&{0}&{0}\\
{0}&a_{22}&{a_{23}}\\
{0}&{0}&a_{33} \end{array} \right).$$ Then, by Eq.(3.9),
$$\left(\begin{array}{cc} a_{22} &a_{23} \\ 0 & a_{33}
\end{array}\right)\in{\mathcal Z}({\mathcal T}_1^\prime),$$
which implies that $$a_{23}=0 \quad{\rm  and}\quad
a_{22}m_{23}=m_{23}a_{33} \eqno(3.10)$$ holds for all
$m_{23}\in{\mathcal M}_{23}$ as ${\mathcal T}_1^\prime$ is a
triangular ring.

Now, with
$$m'=\left(
\begin{array}{ccc}
 0&{m_{12}}&{m_{13}}\\
{0}&0&0\\
{0}&{0}&0 \end{array} \right),$$
by Eq.(3.8), we have $$\begin{array}{rl}0=&[\varphi_{1}(a+m)-\varphi_{1}(a)-\varphi_{1}(m),m']\\
 \if false =&\left(
\begin{array}{ccc}
 a_{11}&{0}&{0}\\
{0}&a_{22}&{a_{23}}\\
{0}&{0}&a_{33} \end{array} \right)\left(
\begin{array}{ccc}
 0&{m_{12}}&{m_{13}}\\
{0}&0&0\\
{0}&0&0 \end{array} \right)\\
&-\left(
\begin{array}{ccc}
 0&{m_{12}}&{m_{13}}\\
{0}&0&0\\
{0}&0&0\end{array} \right)\left(
\begin{array}{ccc}
 a_{11}&{0}&{0}\\
{0}&a_{22}&{a_{23}}\\
{0}&{0}&a_{33} \end{array} \right) \\ \fi =&\left(
\begin{array}{ccc}
 0&{a_{11}m_{12}-m_{12}a_{22}}&{a_{11}m_{13}-m_{13}a_{33}}\\
{0}&0&0\\
{0}&0&0\end{array} \right).\\\end{array}$$ So
$$a_{11}m_{12}-m_{12}a_{22}=0 \quad{\rm and}\quad a_{11}m_{13}-m_{13}a_{33}
=0 \eqno(3.11)$$ hold for all $m_{12}\in{\mathcal M}_{12}$ and all
$m_{13}\in{\mathcal M}_{13}$. This, together with  Eq.(3.10) and
Eq.(3.11) entails that
$$\varphi_{1}(a+m)-\varphi_{1}(a)-\varphi_{1}(m)\in {\mathcal
Z}(\mathcal T). $$

\if false Thus, $a_{11}\in Z(\mathcal{T})$ and
$$[\varphi_{1}(a+m')-\varphi_{1}(a)-\varphi_{1}(m'),a']=0
\eqno(13)$$ for all $a, a'\in\mathcal{T}_{1}, b\in\mathcal{T}'_{1},
m, m'\in M_{1}$. According to (9),(12),(13), we get
$$[\varphi_{1}(a+m')-\varphi_{1}(a)-\varphi_{1}(m'),\mathcal{T}]=0.$$
So
$$\varphi_{1}(a+m')-\varphi_{1}(a)-\varphi_{1}(m')\in Z(\mathcal{T})$$
for all $a, a'\in\mathcal{T}_{1}, b\in\mathcal{T}'_{1}, m, m'\in
M_{1}$. \fi

We remark here that, for $b\in{\mathcal T}_1^\prime$ and
$m\in{\mathcal M}_1$, though a similar argument by using Lemmas
3.3-3.5 as above gives
$$ [\varphi_1(b+m)-\varphi_1(b)-\varphi_1(m), m']=0 \quad \forall
m'\in {\mathcal M}_1
$$
and
$$ [\varphi_1(b+m)-\varphi_1(b)-\varphi_1(m), a']=0 \quad \forall
a'\in {\mathcal T}_1,
$$
\if false Eq.(3.13) means that we can write
$$\varphi_1(b+m)-\varphi_1(b)-\varphi_1(m)=\left(
\begin{array}{ccc}
 b_{11}&{0}&{0}\\
{0}&b_{22}&{b_{23}}\\
{0}&{0}&b_{33} \end{array} \right),
$$
Then Eq.(3.11) implies that
$$ b_{11}m_{12}=m_{12}b_{22} \quad{\rm and}\quad
b_{11}m_{13}=m_{12}b_{23}+m_{13}b_{33}
$$
hold for all $m_{12}\in{\mathcal M}_{12}$ and $m_{13}\in{\mathcal
M}_{13}$. Taking $m_{12}=0$ gives $b_{11}m_{13}=m_{13}b_{33}$ holds
for all $m_{13}\in{\mathcal M}_{13}$, which implies that
$m_{12}b_{23}=0$ for all $m_{12}\in{\mathcal M}_{12}$.\fi we cannot
drive out $\varphi_{1}(b+m)-\varphi_{1}(b)-\varphi_{1}(m)\in\mathcal
Z(\mathcal{T}) $ at this point.

To show that
$\varphi_{1}(a+b)-\varphi_{1}(a)-\varphi_{1}(b)\in\mathcal
Z(\mathcal{T}) $, we need prove (2) first.

(2)  By Lemma 3.4, $\varphi_{1}(Q_{1}),\varphi_{1}(Q_{1}')\in
{\mathcal Z}(\mathcal{T})$. As we just have proved in (1), we have
$\varphi_{1}(Q_{1}+m)-\varphi_{1}(Q_{1})-\varphi_{1}(m)\in \mathcal
Z(\mathcal{T})$. Therefore,
$$[\varphi_{1}(Q_{1}+m),m'+Q'_{1}]=[\varphi_{1}(m),m'+Q'_{1}]=\varphi_{1}(m)\eqno(3.12)$$ and
$$[Q_{1}+m,\varphi_{1}(m'+Q'_{1})]=[Q_{1}+m,\varphi_{1}(m')]=\varphi_{1}(m')\eqno(3.13)$$
hold for all $ m, m'\in \mathcal M_{1}$.

 Since
$m+m'=[Q_{1}+m,m'+Q'_{1}]$,  it follows from Eqs.(3.12) and (3.13)
that
$$\begin{array}{rl}\varphi_{1}(m+m')=&\varphi_{1}([Q_{1}+m, m'+Q'_{1}])\\
=&[\varphi_{1}(Q_{1}+m),m'+Q'_{1}]+[Q_{1}+m,\varphi_{1}(m'+Q'_{1})]\\
=&\varphi_{1}(m)+\varphi_{1}(m')  \end{array}$$ for all $ m, m'\in
\mathcal M_{1}$. Hence, $\varphi_{1}$ is additive on $\mathcal
M_{1}$.

Now let us go back to the proof of
$\varphi_{1}(a+b)-\varphi_{1}(a)-\varphi_{1}(b)\in\mathcal
Z(\mathcal{T}) $.

  Let $a\in{\mathcal T}_1$, $b\in{\mathcal T}_1^\prime$. For any
$m\in{\mathcal M}_1$, as $[a+b,m]=am-mb$, by (2), we have
$$\begin{array}{rl}&[\varphi_{1}(a+b),m]+[a+b,\varphi_{1}(m)]=\varphi_{1}([a+b,m])=\varphi_{1}(am-mb)\\
= &\varphi_{1}(am)-\varphi_{1}(mb)= \varphi_{1}([a,m])+\varphi_{1}([b,m])\\
=&[\varphi_{1}(a),m]+[a,\varphi_{1}(m)]+[\varphi_{1}(b),m]+[b,\varphi_{1}(m)].
\end{array}$$
So
$$[\varphi_{1}(a+b)-\varphi_{1}(a)-\varphi_{1}(b),m]=0 \eqno(3.14)$$
holds for all $m\in{\mathcal M}_1$. Also, for any
$b^\prime\in{\mathcal T}_1^\prime$,
$$\begin{array}{rl}& [\varphi_{1}(a+b),b']+[a+b,\varphi_{1}(b')]=\varphi_{1}([a+b,b'])\\
=&\varphi_{1}([b,b'])
=[\varphi_{1}(b),b])]+[b,\varphi_{1}(b')].\end{array}$$ This entails
that
$$[\varphi_{1}(a+b)-\varphi_{1}(a)-\varphi_{1}(b),b']=0 \eqno(3.15)$$
holds for all $b^\prime \in{\mathcal T}_1^\prime$. Similarly one can
check that
$$[\varphi_{1}(a+b)-\varphi_{1}(a)-\varphi_{1}(b),a']=0 \eqno(3.16)$$
holds for all $a^\prime \in{\mathcal T}_1$. Eqs.(3.14)-(3.16)
together imply that
$$[\varphi_{1}(a+b)-\varphi_{1}(a)-\varphi_{1}(b),\mathcal{T}]=0,$$
and hence, $$\varphi_{1}(a+b)-\varphi_{1}(a)-\varphi_{1}(b)\in
{\mathcal Z}(\mathcal{T}) \eqno(3.17)$$ holds for all $a\in{\mathcal
T}_1$ and $b\in{\mathcal T}_1^\prime$.

(3) Now consider the general case of $a+m+b$. Note that
$Q_{1}\varphi_{1}(a+b+m)Q'_{1}=\varphi_{1}(m)$. For any
$a'\in\mathcal T_1$, since
$$\varphi_{1}([a+b+m,a'])=[\varphi_{1}(a+b+m),a']+[a+b+m,\varphi_{1}(a')]$$
and
$$\begin{array}{rl}\varphi_{1}([a+b+m,a'])&=\varphi_{1}([a+m,a'])\\
&=[\varphi_{1}(a+m),a']+[a+m,\varphi_{1}(a')]. \end{array}$$   we
get
$$[\varphi_{1}(a+b+m)-\varphi_{1}(a+m),a']=0 \quad \forall a'\in{\mathcal T}_1.$$
This, together with
$\varphi_{1}(a+m)-\varphi_{1}(a)-\varphi_{1}(m)\in{\mathcal
Z}(\mathcal T)$  by (1) and an obvious fact $[\varphi(b),a']=0$ by
Lemma 3.3, gives
$$[\varphi_{1}(a+b+m)-\varphi_{1}(a)-\varphi_{1}(b)-\varphi_{1}(m),a']=0, \quad \forall a'\in\mathcal T_1. \eqno(3.18)$$
Similarly, for any $m'\in{\mathcal M}_1$,
$$\varphi_{1}([a+b+m,m'])=[\varphi_{1}(a+b+m),m']+[a+b+m,\varphi_{1}(m')]$$
and
$$\begin{array}{rl}\varphi_{1}([a+b+m,m'])&=\varphi_{1}([a+b,m'])\\
&=[\varphi_{1}(a+b),m']+[a+b,\varphi_{1}(m')] \end{array}$$ leads to
$$[\varphi_{1}(a+b+m)-\varphi_{1}(a+b),m']=0.$$
Then by Eq.(3.17) proved before we obtain that
$$[\varphi_{1}(a+b+m)-\varphi_{1}(a)-\varphi_{1}(b)-\varphi_{1}(m),m']=0 \quad \forall m'\in{\mathcal M}_1.\eqno(3.19)$$
Write $\varphi_{1}(a+b+m)$, $\varphi_{1}(a)$, $\varphi_{1}(b)$,
$\varphi_{1}(m)$ and $m'$ in the form
$$\varphi_{1}(a+b+m)=\left(
\begin{array}{ccc}
 c_{11}&c_{12}&c_{13}\\
{0}&c_{22}&{c_{23}}\\
{0}&{0}&c_{33} \end{array} \right),\quad \varphi_{1}(b)=\left(
\begin{array}{ccc}
 b_{11}&0&0\\
{0}&b_{22}&{b_{23}}\\
{0}&{0}&b_{33} \end{array} \right),$$ $$ \varphi_{1}(a)=\left(
\begin{array}{ccc}
 a_{11}&0&0\\
{0}&a_{22}&{a_{23}}\\
{0}&{0}&a_{33} \end{array} \right),\quad \varphi_{1}(m)=\left(
\begin{array}{ccc}
 0&m_{12}&{m_{13}}\\
{0}&0&0\\
{0}&{0}&0 \end{array} \right)$$ and $$ m'=\left(
\begin{array}{ccc}
 0&m'_{12}&{m'_{13}}\\
{0}&0&0\\
{0}&{0}&0 \end{array} \right).$$
Then
$$\begin{array}{rl}& [\varphi_{1}(a+b+m),m'] \\ =& \left(
\begin{array}{ccc}
 c_{11}&c_{12}&c_{13}\\
{0}&c_{22}&{c_{23}}\\
{0}&{0}&c_{33} \end{array} \right)\left(
\begin{array}{ccc}
 0&m'_{12}&{m'_{13}}\\
{0}&0&0\\
{0}&{0}&0 \end{array} \right)-\left(
\begin{array}{ccc}
 0&m'_{12}&{m'_{13}}\\
{0}&0&0\\
{0}&{0}&0 \end{array} \right)\left(
\begin{array}{ccc}
 c_{11}&c_{12}&c_{13}\\
{0}&c_{22}&{c_{23}}\\
{0}&{0}&c_{33} \end{array} \right)\\
=& \left(
\begin{array}{ccc}
 0&{c_{11}m'_{12}-m'_{12}a_{22}}&{c_{11}m'_{13}-m'_{13}c_{33}-m'_{12}c_{23}}\\
{0}&0&0\\
{0}&{0}&0 \end{array} \right),\\\end{array}$$
$$\begin{array}{rl} &[\varphi_{1}(b),m']\\= &\left(
\begin{array}{ccc}
 b_{11}&0&0\\
{0}&b_{22}&{b_{23}}\\
{0}&{0}&b_{33} \end{array} \right)\left(
\begin{array}{ccc}
 0&m'_{12}&{m'_{13}}\\
{0}&0&0\\
{0}&{0}&0 \end{array} \right)-\left(
\begin{array}{ccc}
 0&m'_{12}&{m'_{13}}\\
{0}&0&0\\
{0}&{0}&0 \end{array} \right)\left(
\begin{array}{ccc}
 b_{11}&0&0\\
{0}&b_{22}&{b_{23}}\\
{0}&{0}&b_{33} \end{array} \right)\\
=&\left(
\begin{array}{ccc}
 0&{b_{11}m'_{12}-m'_{12}b_{22}}&{b_{11}m'_{13}-m'_{13}b_{33}-m'_{12}b_{23}}\\
{0}&0&0\\
{0}&{0}&0 \end{array} \right)\\\end{array}$$ and
$$\begin{array}{rl} &[\varphi_{1}(a),m']\\= &\left(
\begin{array}{ccc}
 a_{11}&0&0\\
{0}&a_{22}&{a_{23}}\\
{0}&{0}&a_{33} \end{array} \right)\left(
\begin{array}{ccc}
 0&m'_{12}&{m'_{13}}\\
{0}&0&0\\
{0}&{0}&0 \end{array} \right)-\left(
\begin{array}{ccc}
 0&m'_{12}&{m'_{13}}\\
{0}&0&0\\
{0}&{0}&0 \end{array} \right)\left(
\begin{array}{ccc}
 a_{11}&0&0\\
{0}&a_{22}&{a_{23}}\\
{0}&{0}&a_{33} \end{array} \right)\\
=&\left(
\begin{array}{ccc}
 0&{a_{11}m'_{12}-m'_{12}a_{22}}&{a_{11}m'_{13}-m'_{13}a_{33}-m'_{12}a_{23}}\\
{0}&0&0\\
{0}&{0}&0 \end{array} \right).\\\end{array}$$ It follows that
$$(c_{11}-a_{11}-b_{11})m'_{12}-m'_{12}(c_{22}-a_{22}-b_{22})=0$$ and
$$(c_{11}-a_{11}-b_{11})m'_{13}-m'_{13}(c_{33}-a_{33}-b_{33})-m'_{12}(c_{23}-a_{23}-b_{23})=0$$
hold for all $m'_{12}\in{\mathcal M}_{12}$ and $ m'_{13}\in
{\mathcal M}_{13}$. Consequently,
$$(c_{11}-a_{11}-b_{11})m'_{12}=m'_{12}(c_{22}-a_{22}-b_{22}),$$
$$(c_{11}-a_{11}-b_{11})m'_{13}=m'_{13}(c_{33}-a_{33}-b_{33})$$ and
$$m_{12}'(c_{23}-a_{23}-b_{23})=0 $$ hold for all
$m_{12}',m_{13}'\in{\mathcal M}_{12}$. Particularly, we have
$$(c_{11}-a_{11}-b_{11})\in{\mathcal Z}(\mathcal T_{1}).$$ Also, by
Eq.(3.18) we see that $$m_{12}=c_{12}\quad{\rm and}\quad
m_{13}=c_{13}. $$ Let ${\mathcal V}_{23}=\{ w_{23}\in{\mathcal
M}_{23} : {\mathcal M}_{12}w_{23}=0\}.$ Therefore,
$$\begin{array}{rl} & \varphi_{1}(a+m+b)-\varphi_{1}(a)-\varphi_{1}(b)-\varphi_{1}(m) \\ \in
& \{ \left(
\begin{array}{ccc}
 a_{11}&0&0\\
{0}&a_{22}&w_{23}\\
{0}&{0}&a_{33} \end{array} \right) : w_{23}\in{\mathcal V}_{23},
a_{11}\in{\mathcal Z}({\mathcal T}_{1}), a_{11}m_{1j}=m_{1j}a_{jj}
\forall m_{1j}\in\mathcal M_{1j}, j=2,3\}.\end{array}$$ $\hfill\Box$

We remark that,   Lemmas 3.2-3.6 hold for all multiplicative Lie
derivations $\varphi_1$ on ${\mathcal T}$ satisfying
$Q_1\varphi(Q_1')Q_1'=0$, and the assumption  ``$Q_{1}{\mathcal
Z}(\mathcal{T})Q_{1}={\mathcal Z}(\mathcal{T}_{1})$ and
$Q'_{1}{\mathcal Z}(\mathcal{T})Q'_{1}={\mathcal
Z}(\mathcal{T}'_{1})$" as well as the assumption  ``every
multiplicative Lie derivation on ${\mathcal T}_1'$ has the standard
form" are not needed.

 {\bf Proof of Theorem 3.1: $k=1$.}  Note that, by Lemma 3.2 and Lemma 3.3, we have
$Q_{1}\varphi_{1}(Q'_{1})Q'_{1}=0$,
$Q_{1}\varphi_{1}(\mathcal{T}'_{1})Q_{1}\subseteq Q_{1}{\mathcal
Z}(\mathcal{T}) Q_{1}$,
$Q'_{1}\varphi_{1}(\mathcal{T}_{1})Q'_{1}\subseteq Q'_{1}{\mathcal
Z}(\mathcal{T}) Q'_{1}$ and $\varphi_{1}(Q'_{1}),
\varphi_{1}(Q_{1})\in {\mathcal Z}(\mathcal{T})$.

Let $\tau:\pi_{\mathcal{T}_{1}}(\mathcal Z(\mathcal{T}))\rightarrow
\pi_{\mathcal{T}'_{1}}(\mathcal Z(\mathcal{T}))$ be the  unique ring
isomorphism so that $z\oplus \tau(z)\in{\mathcal Z}(\mathcal T)$.
According to the hypotheses on $\mathcal T$ in Theorem 3.1, we have
$\pi_{\mathcal{T}_{1}}(\mathcal Z(\mathcal{T}))={\mathcal
Z}(\mathcal{T}_{1})$ and $\pi_{\mathcal{T}'_{1}}(\mathcal
Z(\mathcal{T}))={\mathcal Z}(\mathcal{T}'_{1})$. Thus, for each
$a\in \mathcal Z(\mathcal{T}_1)$, $am=m\tau(a)$ holds for all $m\in
\mathcal M_{1}$. Clearly, by identifying $\pi_{\mathcal
 T_1}(\mathcal T)$ with $Q_1{\mathcal T}Q_1$ and $\pi_{\mathcal
 T_1^\prime}(\mathcal T)$ with $Q_1'{\mathcal T}Q_1'$, we have
$Q_{1}\varphi_{1}(b)Q_{1}+ \tau (Q_{1}\varphi_{1}(b)Q_{1})\in
{\mathcal Z}(\mathcal{T})$ and
$\tau^{-1}(Q'_{1}\varphi_{1}(a)Q'_{1})+
Q'_{1}\varphi_{1}(a)Q'_{1}\in {\mathcal Z}(\mathcal{T})$ for all
$a\in\mathcal{T}_{1}$ and $ b\in\mathcal{T}'_{1}.$ Let us define
maps $\gamma_{1}:\mathcal{T}\rightarrow {\mathcal Z}(\mathcal{T})$
and $\delta_{0}:\mathcal{T}\rightarrow\mathcal{T}$ by
$$\begin{array}{rl} \gamma_{1}(x)= &
Q_{1}\varphi_{1}(Q'_{1}xQ'_{1})Q_{1}+ \tau
(Q_{1}\varphi_{1}(Q'_{1}xQ'_{1})Q_{1})
 \\ &+ \tau^{-1}(Q'_{1}\varphi_{1}(Q_{1}xQ_{1})Q'_{1})
  +Q'_{1}\varphi_{1}(Q_{1}xQ_{1})Q'_{1}   \end{array}$$ and
$$\delta_{0}(x)=\varphi_{1}(x)-\gamma_{1}(x).$$
Obviously, $ \gamma_{1}({\mathcal M}_{1})=0$. Hence, according to
Lemma 3.3 (c), we have
$$\delta_{0}(m)=\varphi_{1}(m)\in\mathcal M_{1} $$
for all $m\in {\mathcal M}_{1}.$

We claim that $\gamma_{1}([\mathcal{T},\mathcal{T}])=0$ and that
$\delta_{0}$ is a multiplicative Lie derivation.

In fact, for arbitrary $x_{1}, x_{2}\in \mathcal{T},$
$$\begin{array}{rl}\gamma_{1}([x_{1},x_{2}])=&Q_{1}\varphi_{1}( Q'_{1}([x_{1},x_{2}])Q'_{1})Q_{1}+\tau^{-1}(Q'_{1}\varphi_{1}(Q_{1}([x_{1},x_{2}])Q_{1})Q'_{1})\\
&+Q'_{1}\varphi_{1}(Q_{1}([x_{1},x_{2}])Q_{1})Q'_{1}+\tau (Q_{1}\varphi_{1}(Q'_{1}([x_{1},x_{2}]Q'_{1})Q_{1})\\
=&Q_{1}\varphi_{1}([Q'_{1}x_{1}Q'_{1},Q'_{1}x_{2}Q'_{1}])Q_{1}+\tau^{-1}(Q'_{1}\varphi_{1}([Q_{1}x_{1}Q_{1},Q_{1}x_{2}Q_{1}])Q'_{1})\\
&+Q'_{1}\varphi_{1}([Q_{1}x_{1}Q_{1},Q_{1}x_{2}Q_{1}])Q'_{1}+\tau
(Q_{1}\varphi([Q'_{1}x_{1}Q'_{1},Q'_{1}x_{2}Q'_{1}])Q_{1}).
  \end{array}$$
Since
$$\begin{array}{rl}& Q_{1}\varphi_{1}( [Q'_{1}x_{1}Q'_{1},Q'_{1}x_{2}Q'_{1}])Q_{1}\\
=& Q_{1}([\varphi_{1}(Q'_{1}x_{1}Q'_{1}),Q'_{1}x_{2}Q'_{1}])Q_{1}+Q_{1}([Q'_{1}x_{1}Q'_{1},\varphi_{1}(Q'_{1}x_{2}Q'_{1}]))Q_{1}$$\\
=&0
 \end{array}$$
 and
$$\begin{array}{rl}& Q'_{1}\varphi_{1}( [Q_{1}x_{1}Q_{1},Q_{1}x_{2}Q_{1}])Q'_{1}\\
=&  Q'_{1}([\varphi_{1}(Q_{1}x_{1}Q_{1}),Q_{1}x_{2}Q_{1}])Q'_{1}+Q'_{1}([Q_{1}x_{1}Q_{1},\varphi_{1}(Q_{1}x_{2}Q_{1}]))Q'_{1}$$\\
=&0,
 \end{array}$$
we see that  $\gamma_{1}([x_{1},x_{2}])=0$. So
$\gamma_{1}([\mathcal{T},\mathcal{T}])=0.$  Consequently
$$\begin{array}{rl}\delta_{0}([x_{1},x_{2}])&=\varphi_{1}([x_{1},x_{2}])\\
&=[\varphi_{1}(x_{1}),x_{2}]+[x_{1},\varphi_{1}(x_{2})]\\
&=[\varphi_{1}(x_{1})-\gamma_{1}(x_{1}),x_{2}]+[x_{1},\varphi_{1}(x_{2})-\gamma_{1}(x_{2})]\\
&=[\delta_{0}(x_{1}),x_{2}]+[x_{1},\delta_{0}(x_{2})]
 \end{array}$$
for all $x_{1}, x_{2}\in\mathcal{T}$; that is, $\delta_{0}$ is a
multiplicative Lie derivation which obviously satisfies
$Q_{1}\delta_{0}(Q'_{1})Q_{1}=0.$ Furthermore, we have
$$\delta_{0}({\mathcal M}_{1})\subseteq {\mathcal M}_{1},\ \
\delta_{0}(\mathcal{T}_{1})\subseteq \mathcal{T}_{1} \ \ {\rm and}\
\ \delta_{0}(\mathcal{T}'_{1})\subseteq \mathcal{T}'_{1}.
\eqno(3.20)$$ Indeed, $\delta_{0}({\mathcal M}_{1})\subseteq
{\mathcal M}_{1}$ is clear. Using Lemma 3.3 (c) and the definitions
of
  $\delta_{0}$, $\gamma_{1}$, for each $a\in\mathcal{T}$, one gets
$$\begin{array}{rl}\delta_{0}(a)=&\varphi_1(a)-\gamma _1(a) \\
=&Q_{1}\varphi_1(a)(a)Q_{1}+Q'_{1}\varphi_{1}(a)Q'_{1}-\tau^{-1}(Q'_{1}\varphi_{1}(a)Q'_{1})-Q'_{1}\varphi_{1}(a)Q'_{1}\\
=&Q_{1}\varphi_1(a)Q_{1}-\tau^{-1}(Q'_{1}\varphi_{1}(a)Q'_{1})
\in\mathcal{T}_{1}.
 \end{array}$$
 Analogously, one has $\delta_{0}(b)=Q'_{1}\delta_{0}(b)Q'_{1}\in\mathcal{T}'_{1}$ for each $b\in\mathcal{T}'_{1}$.

Next, let us introduce two more maps   $\delta_{1}:
\mathcal{T}\rightarrow\mathcal{T}$ and $\xi_{1}:
\mathcal{T}\rightarrow\mathcal{T}$  defined respectively by
$$\delta_{1}(x)=\delta_{0}(Q_{1}xQ_{1})+\delta_{0}(Q_{1}xQ'_{1})+\delta_{0}(Q'_{1}xQ'_{1})$$
and
$$\xi_{1}(x)=\delta_{0}(x)-\delta_{1}(x)$$
for each $x\in\mathcal{T}.$ Note that $\delta_{0}$ is a
multiplicative Lie derivation. We claim that $\delta_{1}$ is also a
multiplicative Lie derivation. To see this, letting
$a,a'\in\mathcal{T}_{1}, b,b'\in\mathcal{T}'_{1},
m,m'\in\mathcal{M}_{1}$, we have
$$\begin{array}{rl} &\delta_{1}([a+b+m,a'+b'+m'])\\ =& \delta_{1}(\left(
\begin{array}{cc}
a&m\\
{0}&b
 \end{array} \right)\left(
\begin{array}{cc}
a'&m'\\
{0}&b'
 \end{array} \right)-\left(
\begin{array}{cc}
a'&m'\\
{0}&b'
 \end{array} \right)\left(
\begin{array}{cc}
a&m\\
{0}&b
 \end{array} \right))\\
=&\delta_{1}(\left(
\begin{array}{cc}
aa'-a'a&am'+mb'-a'm-m'b\\
{0}&bb'-b'b
 \end{array} \right))\\
=&\delta_{0}(aa'-a'a)+\delta_{0}(am'+mb'-a'm-m'b)+\delta_{0}(bb'-b'b)\\
=&\delta_{0}{([a,a'])}+[\delta_{0}(a),m']+[a,\delta_{0}(m')]+[\delta_{0}(m),b']+[m,\delta_{0}(b')]\\
&+[\delta_{0}(m),a']+[m,\delta_{0}(a')]
+[\delta_{0}(b),m']+[b,\delta_{0}(m')]+\delta_{0}([b,b'])
\end{array}$$
and
$$\begin{array}{rl} & [\delta_{1}(a+b+m),a'+b'+m']
+[a+b+m,\delta_{1}(a'+b'+m')]\\
=&[\delta_{0}(a)+\delta_{0}(b)+\delta_{0}(m),a'+b'+m']+[a+b+m,\delta_{0}(a')+\delta_{0}(b')+\delta_{0}(m')]\\
=&[\delta_{0}(a),a']+[\delta_{0}(a),m']+[\delta_{0}(b),b']+[\delta_{0}(b),m']+[\delta_{0}(m),a']+[\delta_{0}(m),b']\\
&+[a,\delta_{0}(a')]+[m,\delta_{0}(a')]+[b,\delta_{0}(b')]+[m,\delta_{0}(b')]+[a,\delta_{0}(m')]+[b,\delta_{0}(m')].
 \end{array}$$
So $$\delta_{1}([a+b+m,a'+b'+m'])=[\delta_{1}(a+b+m),a'+b'+m']
+[a+b+m,\delta_{1}(a'+b'+m')] .$$
It follows that $\delta_{1}$ is a
multiplicative Lie derivation with
$Q_{1}\delta_{1}(Q'_{1})Q'_{1}=0.$ Consequently,
$\xi_{1}(x)=\delta_{0}(x)-\delta_{1}(x)$ is also a multiplicative
Lie derivation with $Q_{1}\xi_{1}(Q'_{1})Q'_{1}=0.$ This ensures
that Lemma 3.3-Lemma 3.5 are applicable to both $\delta_1$ and
$\xi_1$.

Our next aim is to prove that $\delta_{1}$ is a derivation of
$\mathcal{T}$.

Let $a,a'\in\mathcal{T}_{1}$. For each $m\in {\mathcal M}_{1}$.
Using Lemmas 3.4-3.6 and the fact that
$\delta_{1}(\mathcal{T}_{1})\subseteq\mathcal{T}_{1}$ (see
Eq.(3.20)), we get
$$\begin{array}{rl}\delta_{1}((a+a')m)&=\delta_{1}(am+a'm)\\
&=\delta_{1}(am)+\delta_{1}(a'm)\\
&=\delta_{1}(a)m+a\delta_{1}(m)-m\delta_{1}(a)+\delta_{1}(a')m+a'\delta_{1}(m)-m\delta_{1}(a')\\
&=\delta_{1}(a)m+a\delta_{1}(m)+\delta_{1}(a')m+a'\delta_{1}(m).
\end{array}
 $$ On the other hand, by using Lemma 3.5 we have
$$\begin{array}{rl}\delta_{1}((a+a')m)&=\delta_{1}(a+a')m+(a+a')\delta_{1}(m)-m\delta_{1}(a+a')\\
&=\delta_{1}(a+a')m+a\delta_{1}(m)+a'\delta_{1}(m).
\end{array}
 $$
The above two equations ensure that
$$(\delta_{1}(a+a')-\delta_{1}(a)-\delta_{1}(a'))m=0 \eqno(3.21)$$
holds for all $a, a'\in\mathcal{T}_{1}$ and $m\in{\mathcal M}_1$. By
the definition of triangular 3-matrix ring, it is clear that
${\mathcal M}_{1}$ is faithful as a left $\mathcal{T}_{1}$-module.
So Eq.(3.21) entails that $
\delta_{1}(a+a')=\delta_{1}(a)+\delta_{1}(a')$ holds for all $a,
a'\in\mathcal{T}_{1}$; that is, $\delta_{1}$ is additive on
$\mathcal{T}_{1}$.

In addition, using Lemma 3.5 on $\delta_{1}(aa'm)$, we see that
$$\begin{array}{rl}\delta_{1}(aa'm)&=\delta_{1}(aa')m+aa'\delta_{1}(m)-m\delta_{1}(aa')\\
&=\delta_{1}(aa')m+aa'\delta_{1}(m).
\end{array}
$$
While on the other hand
$$\begin{array}{rl}\delta_{1}(aa'm)=& \delta_{1}(a(a'm)) =\delta_{1}(a)a'm+a\delta_{1}(a'm)-a'm\delta_{1}(a)\\
=&\delta_{1}(a)a'm+a\delta_{1}(a')m+aa'\delta_{1}(m)-am\delta_{1}(a')-a'm\delta_{1}(a)\\
=& \delta_{1}(a)a'm+a\delta_{1}(a')m+aa'\delta_{1}(m).
\end{array}
$$
These force that
 $\delta_{1}(aa')=\delta_{1}(a)a'+a\delta_{1}(a')$ and hence $\delta_{1}$ is a  derivation when restricted to $\mathcal{T}_{1}$.

In the sequel we denote  $\delta_{1}|_{\mathcal{T}_{1}}$ by
$\delta_{\mathcal{T}_{1}}$ and $\delta_{1}|_{\mathcal M_{1}}$ by
$\delta_{M_{1}}$.

Note that $\mathcal{T}'_{1}$ is a triangular ring. By the hypotheses
of Theorem 3.1, the multiplicative Lie derivation
$\delta_{1}|_{\mathcal{T}'_{1}}:
\mathcal{T}'_{1}\rightarrow\mathcal{T}'_{1}$ has a standard form
$\delta_{1}|_{\mathcal{T}'_{1}}=\delta_{\mathcal{T}'_{1}}+\eta$,
where $\delta_{\mathcal{T}'_{1}}$ is an additive derivation on
$\mathcal{T}'_{1}$ and $\eta :\mathcal{T}'_{1}\rightarrow {\mathcal
Z}(\mathcal{T}'_{1})$ is a map such that
$\eta(\mathcal{T}'_{1},\mathcal{T}'_{1})=0$. We assert that
$\eta=0$. In fact, for any $x=\left(
\begin{array}{cc}
 a_{1}&m\\
{0}&a_{2}
\end{array} \right), y=\left(
\begin{array}{cc}
 b_{1}&n\\
{0}&b_{2}
\end{array} \right)\in\mathcal T$, where $a_{1}, b_{1}\in\mathcal{T}_{1}$, $ m, n\in\ {\mathcal M}_{1}$, $ a_{2}, b_{2}\in\mathcal{T}'_{1}$,
 we have
$$\begin{array}{rl}\delta_{1}([x.y])=&\delta_{1}(\left(
\begin{array}{cc}
 a_{1}b_{1}-b_{1}a_{1}&a_{1}n+mb_{2}-b_{1}m-na_{2}\\
{0}&a_{2}b_{2}-b_{2}a_{2}
\end{array} \right))\\
=&\delta_{\mathcal{T}_{1}}( a_{1}b_{1}-b_{1}a_{1})+\delta_{M_{1}}(a_{1}n+mb_{2}-b_{1}m-na_{2})\\
&+\delta_{\mathcal{T}'_{1}}( a_{2}b_{2}-b_{2}a_{2})+\eta( a_{2}b_{2}-b_{2}a_{2}).\\
\end{array}\eqno(3.22)$$
On the other hand
{\small $$\begin{array}{rl} &[\delta_{1}(x),y]+[x,\delta_{1}(y)]=\delta_{1}(x)y-y\delta_{1}(x)+x\delta_{1}(y)-\delta_{1}(y)x\\
=&\left(
\begin{array}{cc}
 \delta_{\mathcal{T}_{1}}(a_{1})& \delta_{M_{1}}(m)\\
{0}& \delta_{\mathcal{T}'_{1}}(a_{2})+\eta(a_{2})
\end{array} \right)\left(
\begin{array}{cc}
 b_{1}&n\\
{0}&b_{2}
\end{array} \right)-\left(
\begin{array}{cc}
 b_{1}&n\\
{0}&b_{2}
\end{array} \right)\left(
\begin{array}{cc}
 \delta_{\mathcal{T}_{1}}(a_{1})& \delta_{M_{1}}(m)\\
{0}& \delta_{\mathcal{T}'_{1}}(a_{2})+\eta(a_{2})
\end{array} \right)\\
&+\left(
\begin{array}{cc}
 a_{1}&m\\
{0}&a_{2}
\end{array} \right)\left(
\begin{array}{cc}
 \delta_{\mathcal{T}_{1}}(b_{1})& \delta_{M_{1}}(n)\\
{0}& \delta_{\mathcal{T}'_{1}}(b_{2})+\eta(b_{2})
\end{array} \right)-\left(
\begin{array}{cc}
 \delta_{\mathcal{T}_{1}}(b_{1})& \delta_{M_{1}}(n)\\
{0}& \delta_{\mathcal{T}'_{1}}(b_{2})+\eta(b_{2})
\end{array} \right)\left(
\begin{array}{cc}
 a_{1}&m\\
{0}&a_{2}
\end{array} \right)\\
=&\left(
\begin{array}{cc}
\delta_{\mathcal{T}_{1}}( a_{1})b_{1}-b_{1}\delta_{\mathcal{T}_{1}}( a_{1})&\delta_{\mathcal{T}_{1}}( a_{1})n+ \delta_{M_{1}}(m)b_{2}-b_{1}\delta_{M_{1}}(m)-n \delta_{\mathcal{T}'_{1}}(a_{2})-n\eta(a_{2})\\
{0}&\delta_{\mathcal{T}'_{1}}(a_{2})b_{2}+\eta(a_{2})b_{2}-b_{2}\delta_{\mathcal{T}'_{1}}(a_{2})-b_{2}\eta(a_{2})
\end{array} \right)\\
&+\left(
\begin{array}{cc}
 a_{1}\delta_{\mathcal{T}_{1}}( b_{1})-\delta_{\mathcal{T}_{1}}( b_{1})a_{1}&a_{1} \delta_{M_{1}}(n)+m\delta_{\mathcal{T}'_{1}}(b_{2})+m\eta(b_{2})-\delta_{\mathcal{T}_{1}}( b_{1})m-\delta_{M_{1}}(n)a_{2}\\
{0}&a_{2}\delta_{\mathcal{T}'_{1}}(b_{2})+a_{2}\eta(b_{2})-\delta_{\mathcal{T}'_{1}}(b_{2})a_{2}-a_{2}\eta(b_{2})\\
\end{array} \right).
\end{array}
\eqno(3.23)$$}
As
$\delta_{1}([x.y])=[\delta_{1}(x),y]+[x,\delta_{1}(y)]$, by
Eq(3.22), Eq(3.23) and   Lemmas 3.4-3.5, we obtain that
$$\begin{array}{rl} &\delta_{\mathcal{T}_{1}}(a_{1})n+a_{1}\delta_{M_{1}}(n)-n\delta_{\mathcal{T}_{1}}(a_{1})+\delta_{M_{1}}(m)b_{2}+m\delta_{\mathcal{T}'_{1}}(b_{2})-\delta_{\mathcal{T}'_{1}}(b_{2})m\\
&-\delta_{\mathcal{T}_{1}}(b_{1})m-b_{1}\delta_{M_{1}}(m)+m\delta_{\mathcal{T}_{1}}(b_{1})-\delta_{M_{1}}(n)a_{2}-n\delta_{\mathcal{T}'_{1}}(a_{2})+\delta_{\mathcal{T}'_{1}}(a_{2})n\\
=&\delta_{M_{1}}(a_{1}n)+\delta_{M_{1}}(mb_{2})-\delta_{M_{1}}(b_{1}m)-\delta_{M_{1}}(na_{2})\\
=&\delta_{M_{1}}(a_{1}n+mb_{2}-b_{1}m-na_{2}) \\
=&\delta_{\mathcal{T}_{1}}( a_{1})n+ \delta_{M_{1}}(m)b_{2}-b_{1}\delta_{M_{1}}(m)-n \delta_{\mathcal{T}'_{1}}(a_{2})-n\eta(a_{2})\\
&+a_{1} \delta_{M_{1}}(n)+m\delta_{\mathcal{T}'_{1}}(b_{2})+m\eta(b_{2})-\delta_{\mathcal{T}_{1}}( b_{1})m-\delta_{M_{1}}(n)a_{2}.\\
\end{array}$$
Since
$n\delta_{\mathcal{T}_{1}}(a_{1})=m\delta_{\mathcal{T}_{1}}(b_{1})=\delta_{\mathcal{T}'_{1}}(a_{2})n=\delta_{\mathcal{T}'_{1}}(b_{2})m=0
$, we see that
$$m\eta(b_{2})-n\eta(a_{2})=0 \quad\mbox{\rm for all}\ n,
m\in{\mathcal M}_{1},$$ which implies that $\eta=0$. To see this,
for any $b_2\in{\mathcal T}'_1$,  picking  $n=0$ in the above
equation gives that $m\eta(b_{2})=0$ holds for all $m\in{\mathcal
M}_1$. Write $m= (
 m_{12}, m_{13})$, $\eta(b_{2})=\left(
\begin{array}{cc}
 z_{2}&0\\
{0}&z_{3}
\end{array} \right)\in{\mathcal Z}({\mathcal T}_1')$. Then
$m\eta(b_{2})=( m_{12}z_{2}, m_{13}z_{3}) =(0, 0)$. As
$m_{12}\in{\mathcal M}_{12}$ and $m_{13}\in{\mathcal M}_{13}$ are
arbitrary, we must have $z_{2}=0$ and $z_{3}=0$. So $\eta(b_{2})=0$.
  Thus $\delta_{1}$ is in fact an
additive derivation on $\mathcal{T}'_{1}$.

By now, we have proved that $\delta_{1}$ is  additive respectively
on $\mathcal{T}_{1}$, $\mathcal{T}'_{1}$, and $\mathcal M_1$ (see
Lemma 3.6 (2)). As $\delta_{1}(\mathcal{T}_{1})\subseteq
\mathcal{T}_{1}, \delta_{1}(\mathcal{T}'_{1})\subseteq
\mathcal{T}'_{1}, \delta_{1}(M_{1})\subseteq M_{1}$, it is easily
checked the $\delta_1$ is additive.
\if false $$\begin{array}{rl}\delta_{1}(x+y)&=\delta_{1}((a_{1}+b_{1})+(m+n)+(a_{2}+b_{2}))\\
&=\delta_{1}(a_{1}+b_{1})+\delta_{1}(m+n)+\delta_{1}(a_{2}+b_{2})\\
&=\delta_{1}(a_{1})+\delta_{1}(b_{1})+\delta_{1}(m)+\delta_{1}(n)+\delta_{1}(a_{2})+\delta_{1}(b_{2})\\
&=\delta_{1}(a_{1}+m+b_{1})+\delta_{1}(a_{2}+n+b_{2})\\
&=\delta_{1}(x)+\delta_{1}(y).\\\end{array}$$ Thus, $\delta_{1}$ is
additive on $\mathcal{T}$.\fi

For any $ x=\left(
\begin{array}{cc}
 a_{1}&m\\
{0}&a_{2}
\end{array} \right)$ and $ y=\left(
\begin{array}{cc}
 b_{1}&n\\
{0}&b_{2}
\end{array} \right)\in\mathcal T$, using Lemmas
3.3-3.5 and the fact  that $\delta_{1}$ is derivation on
$\mathcal{T}_{1}$ and $ \mathcal{T}'_{1}$,    we see that
$$\begin{array}{rl}\delta_{1}(xy)=&\delta_{1}(a_{1}b_{1}+(a_{1}n+mb_{2})+a_{2}b_{2})\\
=&\delta_{1}(a_{1}b_{1})+\delta_{1}(a_{1}n+mb_{2})+\delta_{1}(a_{2}b_{2})\\
=&\delta_{1}(a_{1})b_{1}+a_{1}\delta_{1}(b_{1})+\delta_{1}(a_{1}n)+\delta_{1}(mb_{2})+\delta_{1}(a_{2})b_{2}+a_{2}\delta_{1}(b_{2})\\
=&\delta_{1}(a_{1})b_{1}+a_{1}\delta_{1}(b_{1})+\delta_{1}(a_{1})n+a_{1}\delta_{1}(n)\\
&+\delta_{1}(m)b_{2}+m\delta_{1}(b_{2})+\delta_{1}(a_{2})b_{2}+a_{2}\delta_{1}(b_{2}).\\
\end{array}$$
While on the other hand
$$\begin{array}{rl}\delta_{1}(x)y+x\delta_{1}(y)=&(\delta_{1}(a_{1})+\delta_{1}(m)+\delta_{1}(a_{2}))(b_{1}+n+b_{2})\\
&+(a_{1}+m+a_{2})(\delta_{1}(b_{1})+\delta_{1}(n)+\delta_{1}(b_{2}))\\
=&\delta_{1}(a_{1})b_{1}+\delta_{1}(a_{1})n+\delta_{1}(m)b_{2}+\delta_{1}(a_{2})b_{2}\\
&+a_{1}\delta_{1}(b_{1})+a_{1}\delta_{1}(n)+m\delta_{1}(b_{2})+a_{2}\delta_{1}(b_{2}).\\\end{array}$$
Consequently, $\delta_{1}(xy)=\delta_{1}(x)y+x\delta_{1}(y)$. Hence
$\delta_{1}$ is a  derivation of $\mathcal{T}$.

Finally, let us consider the structure of $\xi_1$.  Since
$\xi_{1}(x)=\delta_{0}(x)-\delta_{1}(x)$ and
$\varphi_{1}(x)=\delta_{0}(x)+\gamma_{1}(x)$, by the definition of
$\gamma_{1}$, we get
$\xi_{1}(x)=\varphi_{1}(x)-\varphi_{1}(Q_{1}xQ_{1})-\varphi_{1}(Q'_{1}xQ'_{1})-\varphi_{1}(Q_{1}xQ'_{1})$.
It follows from Lemma 3.6 (3) that, for each $x\in\mathcal{T},
\xi_{1}(x)\in {\mathcal S}_{23}$, where
$$ \begin{array}{rl} \mathcal S_{23}=& \{   \left(
\begin{array}{ccc}
 a_{11}&0&0\\
{0}&a_{22}&w_{23}\\
{0}&{0}&a_{33} \end{array} \right) : \\  &  w_{23}\in{\mathcal
V}_{23}, a_{11}\in{\mathcal Z}({\mathcal T}_{1}),
a_{11}m_{1j}=m_{1j}a_{jj} \ \forall m_{1j}\in{\mathcal M}_{1j},
j=2,3 \}. \end{array}$$

Define   maps $\gamma_{2}:
\mathcal{T}\rightarrow\mathcal{Z}(\mathcal{T})$ and $\xi_{2}:
\mathcal{T}\rightarrow\mathcal{T}$ by
$$\gamma_{2}(x)=Q_{1}\xi_{1}(x)Q_{1}+\tau(Q_{1}\xi_{1}(x)Q_{1})$$
and
$$\xi_{2} =\xi_{1} -\gamma_{2} .$$
We claim that $\xi_2(x)\in{\mathcal S}_1$ for all $x\in\mathcal T$.
To see this, write $$\xi_1(x)=\left(
\begin{array}{ccc}
 a_{11}&0&0\\
{0}&a_{22}&w_{23}\\
{0}&{0}&a_{33} \end{array} \right)\in{\mathcal S}_{23}\quad{\rm
and}\quad  \gamma_2(x)=\left(
\begin{array}{ccc}
 a_{11}&0&0\\
{0}&b_{22}&0\\
{0}&{0}&b_{33} \end{array} \right)\in{\mathcal Z}(\mathcal T).$$ As
$\xi_1(x)-\gamma_2(x)\in{\mathcal S}_{23}$, we get
$(a_{11}-a_{11})m_{1j}=m_{1j}(a_{jj}-b_{jj})$ holds for all
$m_{1j}\in\mathcal M_{1,j}$, which forces that $b_{jj}=a_{jj}$,
$j=2,3$.

We show that $\gamma_{2}([\mathcal{T},\mathcal{T}])=0$ and hence $\xi_{2}$ is a multiplicative Lie derivation.\\
In fact, for arbitrary $x_{1}, x_{2}\in\mathcal T,$
$$\gamma_{2}([x_{1},x_{2}])=Q_{1}\xi_{1}([x_{1},x_{2}])Q_{1}+\tau(Q_{1}\xi_{1}([x_{1},x_{2}])Q_{1}).$$
Since $Q_1\xi_1(x_i)Q_1\in{\mathcal Z}(\mathcal T_1)$ for $i=1,2$,
we see that
$$\begin{array}{rl}& Q_{1}\xi_{1}([x_{1},x_{2}])Q_{1}\\
=&Q_{1}[\xi_{1}(x_{1}),x_{2}]Q_{1}+Q_{1}[x_{1},\xi_{1}(x_{2})]Q_{1}\\
=& Q_{1}\xi_{1}(x_{1})x_{2}Q_{1}- Q_{1}x_{2}\xi_{1}(x_{1})Q_{1}+Q_{1}x_{1}\xi_{1}(x_{2})Q_{1}- Q_{1}\xi_{1}(x_{2})x_{1}Q_{1}$$\\
=&0,
 \end{array}$$
which entails  $\gamma_{2}([x_{1},x_{2}])=0$. So
$\gamma_{2}([\mathcal T,\mathcal T])=0$. Consequencely
  $\xi_{2}$ is a
multiplicative Lie derivation which obviously satisfies
$Q_{1}\xi_{2}(Q'_{1})Q_{1}'=0$ and $$\xi_{2}\in \{ \left(
\begin{array}{ccc}
 0&0&0\\
{0}&0&w_{23}\\
{0}&{0}&0 \end{array} \right): w_{23}\in{\mathcal V}_{23}
\}={\mathcal S}_1. \eqno(3.24)$$

So far we have shown that $\varphi_1=
\delta_0+\gamma_1=\delta_1+\xi_1+\gamma_1=\delta_1+\xi_2+\gamma_2+\gamma_1$.
Let $\delta=d_1+\delta_1$, $\gamma=\gamma_1+\gamma_2$ and
$\xi=\xi_2$. Then $\delta$ is a derivation, $\gamma$ is a center
valued map vanishing each commutator, $\xi $ is a multiplicative Lie
derivation with range in $\mathcal S_1$ and
$\varphi=\delta+\gamma+\xi$, as desired. \hfill $\Box$

\subsection{Proof of Theorem 3.1: $k=3$.}

In this case $\mathcal T={\mathcal T}_3'+{\mathcal M}_3+\mathcal
T_3$ with ${\mathcal T}_3=Q_3{\mathcal T }Q_3$, ${\mathcal
T}_3'=Q_3'{\mathcal T }Q_3'$ and ${\mathcal M}_3=Q_3'{\mathcal
T}Q_3$.

Let $\varphi:\mathcal T\to\mathcal T$ be a multiplicative Lie
derivation. One can prove in a similar way that the analogues of
Lemmas 3.2-3.5 are valid. We list them here  by omitting the proofs
for application later.

{\bf Lemma 3.2$'$.} {\it There exists an inner derivation
$d_{1}':\mathcal{T}-\mathcal{T}$ and a multiplicative Lie derivation
$\varphi_{1}':\mathcal{T}\rightarrow\mathcal{T}$ such that
$$\varphi=d_{1}'+\varphi_{1}' \quad {\rm and}\quad  Q_{3}'\varphi'_{1}(Q'_{3})Q_{3}=0.$$}

Note that we take $d_1'(x)=[x,\varphi(Q_3')]$.

{\bf Lemma 3.3$'$.} {\it For any $f\in \mathcal{T}'_{3},
g\in\mathcal{T}_{3}$ and $ h\in\mathcal M_{3}$, the following
statements are true:}

(a)   $\varphi'_{1}(Q_{3}),\ \varphi'_{1}(Q'_{3})\in \mathcal
R_{1}\oplus\mathcal R_{2}\oplus\mathcal R_{3}$.

(b)  $Q_{3}'\varphi_{1}(f)Q_{3}=Q'_{3}\varphi'_{1}(g)Q_{3}=0.$

(c)
$\varphi'_{1}(f)=Q'_{3}\varphi'_{1}(f)Q'_{3}+Q_{3}\varphi'_{1}(f)Q_{3}$,
$\varphi'_{1}(g)=Q'_{3}\varphi'_{1}(g)Q'_{3}+Q_{3}\varphi'_{1}(g)Q_{3}$
and $\varphi'_{1}(h)=Q'_{3}\varphi'_{1}(h)Q_{3}.$

(d)  $Q_{3}\varphi'_{1}(\mathcal{T}'_{3})Q_{3}\subseteq \mathcal
Z(\mathcal{T}_{3})$ and $
Q'_{3}\varphi'_{1}(\mathcal{T}_{3})Q'_{3}\subseteq \mathcal
Z(\mathcal{T}'_{3}).$

{\bf Lemma 3.4$'$.} {\it   $\varphi'_{1}(Q_{3}),
\varphi'_{1}(Q'_{3})\in \mathcal Z(\mathcal{T})$ and
$\varphi'_{1}(Q'_{3}xQ_{3})=Q'_{3}\varphi'_{1}(x)Q_{3}$ for all
$x\in\mathcal T$.}

{\bf Lemma 3.5$'$.} {\it For any $f\in \mathcal{T}'_{3}, g\in
\mathcal{T}_{3}$  and $ h\in\mathcal M_{3}$, we have
$$\varphi'_{1}(fh)=\varphi'_{1}(f)h+f\varphi'_{1}(h)-h\varphi'_{1}(f)$$ and $$ \varphi'_{1}(hg)=\varphi'_{1}(h)g+h\varphi'_{1}(g)-\varphi'_{1}(g)h.$$}

 {\bf Lemma 3.6$'$.} {\it For any
$f\in\mathcal{T}_{3}'$, $ g\in\mathcal{T}_{3}$ and $h\in\mathcal
M_{3}$, we have}

(1) {\it Both $\varphi'_{1}(f+h)-\varphi'_{1}(f)-\varphi'_{1}(h)$
and $\varphi'_{1}(f+g)-\varphi'_{1}(f)-\varphi'_{1}(g)\in\mathcal
Z(\mathcal{T}) $.}

(2) {\it $\varphi'_{1}$~is additive on $\mathcal M_{3}$}.

(3) {\it With ${\mathcal V}_{12}=\{ w_{12}\in{\mathcal M}_{12} :
w_{12}{\mathcal M}_{23}=0\}$, we have
$$\begin{array}{rl} &
\varphi'_{1}(f+h+g)-\varphi'_{1}(f)-\varphi'_{1}(g)-\varphi'_{1}(h) \\
\in & \{ \left(
\begin{array}{ccc}
 a_{11}&w_{12}&0\\
{0}&a_{22}&0\\
{0}&{0}&a_{33} \end{array} \right) : w_{12}\in{\mathcal V}_{12},
 a_{ii}m_{i3}=m_{i3}a_{33} \forall m_{i3}\in\mathcal M_{i3}, a_{33}\in{\mathcal Z}({\mathcal T}_{3})\}.\end{array}$$ }

{\bf Proof of Theorem 3.1: $k=3$.} As $Q_{1}{\mathcal
Z}(\mathcal{T})Q_{1}={\mathcal Z}(\mathcal{T}_{1})$ and
$Q'_{1}{\mathcal Z}(\mathcal{T})Q'_{1}={\mathcal
Z}(\mathcal{T}'_{1})$, there is a unique isomorphism $\tau: \mathcal
Z(\mathcal{T}_3)\rightarrow  \mathcal Z(\mathcal{T}_3')$
  so that $\tau(z)+ z \in{\mathcal
Z}(\mathcal T)$.  Then the map $\gamma'_1 $ defined by
$$\begin{array}{rl} \gamma'_{1}(x)= &
Q_{3}\varphi'_{1}(Q'_{3}xQ'_{3})Q_{3}+ \tau'
(Q_{3}\varphi_{1}(Q'_{3}xQ'_{3})Q_{3})
 \\ &+ (\tau')^{-1}(Q'_{3}\varphi'_{1}(Q_{3}xQ_{3})Q'_{3})
  +Q'_{3}\varphi'_{1}(Q_{3}xQ_{3})Q'_{3}   \end{array}$$
  is a center valued map vanishing each commutator. Consequently,
$$\delta'_{0} =\varphi'_{1} -\gamma'_{1}.$$
is still a multiplicative Lie derivation satisfying
$Q_3'\delta_0'(Q'_3)Q_3=0$ and Lemmas 3.2$'$-3.6$'$ are applicable.
By using the hypothesis that every multiplicative Lie derivation on
${\mathcal T}_3'$ has the standard form, one can prove that the map
$\delta_1'$ defined by
$$\delta'_1(f+g+h)= \delta'_0(f)+\delta'_0(g)+\delta'_0(h)$$
for any $f\in{\mathcal T}'_3$, $g\in{\mathcal T}_3$ and
$h\in{\mathcal M}_3$ is a derivation. Thus
$$\xi'_1=\delta_0'-\delta_1'$$ is still a multiplicative Lie
derivation satisfying $Q_3'\xi_1'(Q'_3)Q_3=0$. Moreover, the range
of $\xi'_1$ is contained in the set
$${\mathcal S}_{12}=\{ {\small \left(
\begin{array}{ccc}
 a_{11}&w_{12}&0\\
{0}&a_{22}&0\\
{0}&{0}&a_{33} \end{array} \right) : w_{12}\in{\mathcal V}_{12},
 a_{ii}m_{i3}=m_{i3}a_{33}\ \forall m_{i3}\in\mathcal M_{i3}, i=1,2, a_{33}\in{\mathcal Z}({\mathcal T}_{3})}\}.$$ Let $\gamma'_2$ be
the map defined by
$$\gamma'_2(x)=(\tau')^{-1}(Q_3\xi'_1(x)Q_3)+Q_3\xi'_1(x)Q_3.$$
$\gamma'_2$ is a center valued map vanishing each commutator. Hence
$\xi'_2=\xi'_1-\gamma'_2$ is a multiplicative Lie derivation with
range in ${\mathcal S}_3=\{\left(
\begin{array}{ccc}
0&w_{12}&0\\
{0}&0&0\\
{0}&{0}&0 \end{array} \right) : w_{12}\in{\mathcal V}_{12}\}$ and
$\varphi'_1=\delta'_1+\gamma'_1+\gamma_2'+\xi'_2$. With
$\delta=d_1'+\delta'_1$, $\gamma=\gamma'_1+\gamma'_2$ and
$\xi=\xi'_2$, we see that $\varphi =\delta+\gamma+\xi$ has the
desired form. \hfill$\Box$

\section{Proof of the main result}

We are ready to prove our main result Theorem 2.2.

{\bf Proof of Theorem 2.2}. Assume that $Q_i{\mathcal Z}(\mathcal
T)Q_i={\mathcal Z}(\mathcal T_i)$, $i=1,2,3$.

{\bf Claim 1}. For $k=1,3$, we have $Q_k'{\mathcal Z}(\mathcal
T)Q_k'={\mathcal Z}(\mathcal T'_k)$. Thus $\mathcal T$ satisfies the
hypotheses of Theorem 3.1.

We check $Q_1'{\mathcal Z}(\mathcal T)Q_1'={\mathcal Z}(\mathcal
T'_1)$ and the case $k=3$ is dealt with similarly. Note that
$Q_1'=Q_2+Q_3$ and obviously $Q_1'{\mathcal Z}(\mathcal
J)Q_1'\subseteq {\mathcal Z}(\mathcal T_1')$. Let
$z'=\left(\begin{array}{ccc} 0&0&0\\ 0&z'_{2}&0\\
0&0&z'_3\end{array}\right)\in{\mathcal Z}(\mathcal T'_1)$. Then, for
any $m_{23}\in{\mathcal M}_{23}$ we have $z'_2m_{23}=m_{23}z'_3$.
By the hypothesis of Theorem 2.2, there exists $z=\left(\begin{array}{ccc} z_1&0&0\\ 0&z_{2}&0\\
0&0&z_3\end{array}\right)\in{\mathcal Z}(\mathcal T)$ such that
$z_3=z'_3$. It follows that
$z_2m_{23}=m_{23}z_3=m_{23}z'_3=z'_2m_{23}$ holds for all
$m_{23}\in{\mathcal M}_{23}$ which forces $z_2=z'_2$. So,
$Q_1'zQ_1'=z'$. Hence we have $Q_1'{\mathcal Z}(\mathcal
T)Q_1'={\mathcal Z}(\mathcal T_1')$.

{\bf Claim 2}. If $\xi:\mathcal T\to\mathcal T$ be a multiplicative
Lie derivation of which the range is contained in $\mathcal S_1$,
then $\xi$ is a derivation.

Let $d_{2}: \mathcal{T}\rightarrow\mathcal{T}$ and $\varphi_{2}:
\mathcal{T}\rightarrow\mathcal{T}$  be the maps defined respectively
by
$$d_{2}(x)=[x, \xi(Q'_{3})]$$ and
$$\varphi_{2}(x)=\xi (x)-d_{2}(x)$$
for all $x\in\mathcal{T}$. Obviously, $d_{2}$ is an inner derivation
and  $\varphi_{2}$  is a multiplicative Lie derivation satisfying
$Q_3'\varphi_2(Q_3')Q_3=0$.  Note that $d_2(x)\in {\mathcal S}_1$
and hence $\varphi_2(x)\in {\mathcal S}_1$ for each $x\in\mathcal
T$. Consequently, $\varphi_2(Q'_3)=0$.

We assert that $\varphi_2=0$.

\if false Moreover, for any $f\in\mathcal{T}_{3}',
g\in\mathcal{T}_{3}$ and $ h\in\mathcal{M}_{3}$, we have
$$\varphi_{2}(e+f+n)-\varphi_{2}(e)-\varphi_{2}(f)-\varphi_{2}(n)=\xi_{2}(e+f+n)-\xi_{2}(e)-\xi_{2}(f)-\xi_{2}(n).\eqno(4.1)$$
By Lemma 3.6$'$(3), for any $h'\in\mathcal{M}_{3}$, we have
$$[(\varphi_{2}(f+g+h)-\varphi_{2}(f)-\varphi_{2}(g)-\varphi_{2}(h)),h']=0.$$
So by Eq.(4.1), we have
$$[\xi_{2}(f+g+h)-\xi_{2}(f)-\xi_{2}(g)-\xi_{2}(h),h']=0 \quad
\forall h'\in{\mathcal M}_3.\eqno(4.2)$$  Write
$\xi_{2}(f+g+h)-\xi_{2}(f)-\xi_{2}(g)-\xi_{2}(h)=\left(
\begin{array}{ccc}
 0&0&0\\
{0}&b_{22}&w_{23}\\
{0}&{0}&b_{33} \end{array} \right)$ and $h'=\left(
\begin{array}{ccc}
 0&0&h'_{13}\\
{0}&0&h'_{23}\\
{0}&{0}&0 \end{array} \right)$.
 By Eq.(4.2), we get
$$\begin{array}{rl}0=&[\xi_{2}(f+g+h)-\xi_{2}(f)-\xi_{2}(g)-\xi_{2}(h),h']\\
=&\left(
\begin{array}{ccc}
 0&0&0\\
{0}&b_{22}&w_{23}\\
{0}&{0}&b_{33} \end{array} \right)\left(
\begin{array}{ccc}
 0&0&h'_{13}\\
{0}&0&h'_{23}\\
{0}&{0}&0 \end{array} \right) -\left(
\begin{array}{ccc}
 0&0&h'_{13}\\
{0}&0&h'_{23}\\
{0}&{0}&0 \end{array} \right)\left(
\begin{array}{ccc}
 0&0&0\\
{0}&b_{22}&w_{23}\\
{0}&{0}&b_{33} \end{array} \right)\\
=&\left(
\begin{array}{ccc}
 0&0&h'_{13}b_{33}\\
{0}&0&b_{22}h'_{23}-h'_{23}b_{33}\\
{0}&{0}&0 \end{array} \right). \end{array}$$ It follows that
$h'_{13}b_{33}=0$ and $b_{22}h'_{23}-h'_{23}b_{33}=0$ hold for all
$h'_{13}\in\mathcal M_{13}$ and $h'_{23}\in\mathcal M_{23}$. This
entails that $b_{22}=b_{33}=0$. Thus,
$$ \xi_{2}(f+g+h)-\xi_{2}(f)-\xi_{2}(g)-\xi_{2}(h)   \in \{ \left(
\begin{array}{ccc}
 0&0&0\\
{0}&0&w_{23}\\
{0}&{0}&0 \end{array} \right): w_{23}\in{\mathcal V}_{23}\}. $$
Therefore,
$$\varphi_{2}(e+f+n)-\varphi_{2}(e)-\varphi_{2}(f)-\varphi_{2}(n)\in{\mathcal
V}_{1}. \eqno(4.3)$$

Let  $\delta_{2}: \mathcal{T}\rightarrow\mathcal{T}$ and $\psi:
\mathcal{T}\rightarrow\mathcal{T}$ be defined by
$$\delta_{2}(x)=\varphi_{2}(Q'_{3}xQ'_{3})+\varphi_{2}(Q_{3}xQ_{3})+\varphi_{2}(Q_{3}xQ'_{3})$$
and $$\psi(x)=\varphi_{2}(x)-\delta_{2}(x)$$ for every
$x\in\mathcal{T}$. Notice that $\mathcal T'_3$ as a triangular ring
satisfies the condition $Q_i{\mathcal Z}({\mathcal
T}_3')Q_i={\mathcal Z}(\mathcal T_i)$, $i=1,2$, and thus, by
\cite[Theorem 5.9]{DD}, every multiplicative Lie derivation has the
standard form. As in the   proof of Theorem 3.1 for the case $k=3$,
it is easily checked that $\delta_{2} $ is a derivation and hence
$\psi $ is a multiplicative Lie derivation such that
$Q_3'\psi(Q_3')Q_3=0$ and $\psi(x)\in{\mathcal V}_{1}$ for all
$x\in\mathcal T$. Obviously, $\xi=d_2+\delta_2+\psi=d'+\psi$ with
$d'=d_2+\delta_2$, as desired.\fi

 By Claim 1, we may apply Theorem 3.1 for the case $k=3$  to
 $\varphi_2$. Thus there exists a derivation $\delta'$, a center
 valued map $\gamma'$ and a multiplicative Lie derivation $\xi':
 \mathcal T\to {\mathcal S}_3$ such that
 $\varphi_2=\delta'+\gamma'+\xi'$. However, by checking the proof
 of Theorem 3.1 for the case $k=1$, we must have
 $\delta'(x)\in{\mathcal S}_1$ and $\gamma'=0$. So, for each $x\in\mathcal T$,  $\varphi_2(x)\in
 {\mathcal S}_1\cap{\mathcal S}_3=\{0\}$; that is, $\varphi_2=0$.

 {\bf Claim 3}. Every multiplicative Lie derivation on $\mathcal T$
 has the standard form.

 Let $\varphi :\mathcal T\to\mathcal T$ be a multiplicative Lie derivation. By Claim 1, under the standard  assumption
$Q_i{\mathcal Z}(\mathcal T)Q_i={\mathcal Z}({\mathcal T}_i)$,
$\mathcal T$ meets the assumptions in Theorem 3.1 for both cases
$k=1$ and $k=3$. Then, there exists a derivation, a center valued
map $\gamma$ which vanishes each commutator and a multiplicative Lie
derivation $\xi :\mathcal T\to \mathcal S_1$ such that
$\varphi=\delta +\gamma +\xi$.  However, by Claim 2, $\xi$ is in
fact a derivation. Therefore, $\varphi$ has the standard form.
\hfill$\Box$

\end{document}